\documentclass[12pt]{amsart}

\usepackage[T1]{fontenc} 
\usepackage{amsmath}
\usepackage{amsthm} 
\usepackage{amscd} 
\usepackage{amssymb}
\usepackage{enumerate} 
\usepackage{latexsym} 
\usepackage{xypic}
\usepackage{dsfont} 
\usepackage{array} 
\usepackage{fullpage}

\usepackage{graphicx}
\usepackage{epsfig} 

\newcommand{\chop}{\dag}

\def\epsilon{\varepsilon}

\renewcommand{\phi}{\varphi}

\DeclareMathOperator{\ind}{\mathrm{ind}} 
\DeclareMathOperator{\igeo}{\mathrm{ind_{geo}}}
\DeclareMathOperator{\iq}{\mathrm{ind}_{\CQ}}

\newcommand{\Stab}{\mathrm{Stab}}
\newcommand{\inv}{^{-1}}

\newcommand{\dom}{\mathrm{dom}}
\newcommand{\rank}{\mathrm{rank}} 
\newcommand{\Cayl}{\Gamma}
\newcommand{\conv}{\mathrm{conv}}

\newcommand{\core}{\text{core}} 
\newcommand{\Tt}{\tilde\Gamma}

\newcommand{\FN}{F_N} 

\newcommand{\CVN}{\mathrm{CV}_N}
\newcommand{\barCVN}{\bar{\mathrm{CV}}_N} 
\newcommand{\CQ}{{\mathcal Q}}
\newcommand{\Tobs}{\widehat T^{\mathrm{\scriptstyle obs}}}
\newcommand{\CN}{\mathcal N} 
\newcommand{\CP}{\mathcal P}

\newcommand{\carac}{\mathds{1}}

\newcommand{\CA}{A} 
\newcommand{\CF}{F}
\newcommand{\CK}{S}

\newcommand{\iwip}{iwip}

\newcommand{\R}{\mathbb R}

\newcommand{\N}{\mathbb N}

\def\bar{\overline} 
\def\tilde{\widetilde} 
\def\hat{\widehat}

\newtheorem{thm}{Theorem}[section]
\newtheorem{cor}[thm]{Corollary} 
\newtheorem{lem}[thm]{Lemma}
\newtheorem{prop}[thm]{Proposition}

\newtheorem*{thm*}{Theorem}
\newtheorem*{prop*}{Proposition}
\newtheorem*{thmmainsi}{Theorem~\ref{thm:mainsi}}

\newtheorem*{thmsiiq}{Theorem~\ref{thm:siiq}}
\newtheorem*{thmiq}{Theorem~\ref{thm:iq}}
\newtheorem*{thmsurface}{Theorem~\ref{thm:surface}}

\newtheorem*{thmmaxindex}{Theorem~\ref{thm:maxindex}}

\newtheorem{defn}[thm]{Definition}
\newtheorem*{defn*}{Definition} 
\newtheorem{rem}[thm]{Remark}
\newtheorem*{rem*}{Remark}

\numberwithin{equation}{section} 

\begin{document}
\title{Rips Induction: Index of the dual lamination of an $\R$-tree}

\author{Thierry Coulbois, Arnaud Hilion}

\date{\today }

\begin{abstract}
  Let $T$ be a $\R$-tree in the boundary of the Outer Space $\CVN$,
  with dense orbits. The $\CQ$-index of $T$ is defined by means of the
  dual lamination of $T$. It is a generalisation of the
  Poincar\'e-Lefschetz index of a foliation on a surface. We prove that
  the $\CQ$-index of $T$ is bounded above by $2N-2$, and we study the
  case of equality.  The main tool is to develop the Rips Machine in
  order to deal with systems of isometries on compact $\R$-trees.
  
  Combining our results on the $\CQ$-index with results on the
  classical geometric index of a tree, developed by Gaboriau and
  Levitt \cite{gl-rank}, we obtain a beginning of classification of
  trees.
\end{abstract}

\maketitle

\tableofcontents

\section{Introduction}

The space of minimal, free and discrete actions by isometries of the
free group $\FN$ of finite rank $N\geq 2$ on $\R$-trees has been
introduced by Culler and Vogtmann \cite{cv-moduli}. Its
projectivization is called Outer Space, and we denote it by $\CVN$. It
has a Thurston-boundary $\partial\CVN$, which gives rise to a
compactification $\bar{\mathrm{CV}}_N=\CVN\cup\partial\CVN$. Elements
of this compact space $\bar{\mathrm{CV}}_N$ are projective classes
$[T]$ of minimal, very small actions by isometries of the free group
$\FN$ on $\R$-trees $T$ (see \cite{cl-verysmall} and
\cite{bf-outer}). The reader will find a survey on Outer Space in
\cite{vogt-survey}.

In this article, we introduce and study the $\mathcal{Q}$-index $\iq(T)$ of
$\R$-trees $T$ in $\partial\CVN$ with dense orbits.  The
$\mathcal{Q}$-index of an $\R$-tree (see Section~\ref{subsec:Q-index})
naturally extends the Poincar\'e-Lefschetz index of a foliation on a
surface as explained below.  The main result of our paper regarding
this $\mathcal{Q}$-index is:

\begin{thmiq} 
Let $T$ be an $\R$-tree with a very small, minimal
action of $\FN$ by isometries with dense orbits. 
Then 
\[ 
\iq(T)\leq 2N-2. 
\] 
\end{thmiq}

We also characterize  the case of equality, see Section~\ref{sec:simaxindex}.

This $\CQ$-index charaterizes dynamical properties of trees. Using it
together with the geometric index introduced by Gaboriau and
Levitt~\cite{gl-rank} we obtain a classification of trees. 

\smallskip

Theorem~\ref{thm:iq} already has several important consequences. 

First, it answers a question of Levitt and
Lustig~\cite[Remark~3.6]{ll-north-south} on the finiteness of the
fibres of the map $\CQ$ (see below).

In our paper \cite{ch-b}, we obtain a qualitative classification of
fully irreducible outer automorphisms of free groups which extends
that of Handel and Mosher~\cite{hm-parageometric} and of
Guirardel~\cite{guir-core}. The key point is to interpret the index of
an \iwip\ automorphism \cite{gjll} as the $\CQ$-index of its repelling
tree in $\partial\CVN$.

In our paper with P.~Reynolds~\cite{chr} we define an induction
analogous to Rauzy-Veech for trees in $\partial\CVN$. As we are
working with systems of isometries on compact trees,
Theorem~\ref{thm:iq} is crucially used to ensure that there are points
where to start the splitting procedure.

\subsection{Measured foliations on surfaces}

\begin{figure}[h] 
\begin{center}
\includegraphics{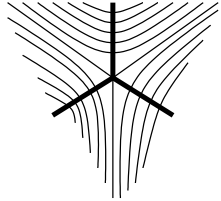} 
\end{center} 
\caption{\label{fig:singularity} 3-prong singularity and the
  transverse tree (in bold).  The local Poincar\'e-Lefschetz index is
  $-\frac 12$ and the local contribution to the $\CQ$-index is $1$.}
\end{figure}

Let $\Sigma$ be a surface of negative Euler characteristic,
$\chi(\Sigma)<0$, with a measured foliation $\mathcal F$ (see
\cite{flp}). The foliation lifts to a measured foliation
$\tilde{\mathcal F}$ of the universal cover $\tilde\Sigma$ of
$\Sigma$. The space of leaves of $\tilde{\mathcal F}$ is an
$\R$-tree $T$: the distance in the tree $T$ is given by the
transverse measure of the foliation $\tilde{\mathcal F}$
(see for instance \cite[chapter 11]{kapo-book}).
This tree comes with a small action of the fundamental
group of $\Sigma$. When $\Sigma$ has non empty
boundary, its fundamental group is a free group $\FN$
and $T$ defines an element of $\bar{\mathrm{CV}}_N$.
The foliation $\mathcal F$ has $k$-prong singularities which give
rise to branch points of valence $k$ in the tree. Locally the
picture is as in Figure~\ref{fig:singularity}.
We say that the foliation and the tree are dual to each other.

A local index can be defined for each singularity $P$: the
Poincar\'e-Lefschetz index. In this paper, we rather consider
minus two times this index: $\ind(P)$ can be defined as the number of 
half-leaves reaching the singular point, minus two. 
Alternatively, $\ind(P)$ is the valence
of the corresponding point in the tree, minus two. Adding-up over
all singular points in $\Sigma$, we obtain the (global) index of
the foliation, which turns out to be equal to $-2\chi(\Sigma)$
($=2N-2$ when $\pi_1(\Sigma)=\FN$).

\begin{figure}[h] 
\begin{center}
\includegraphics{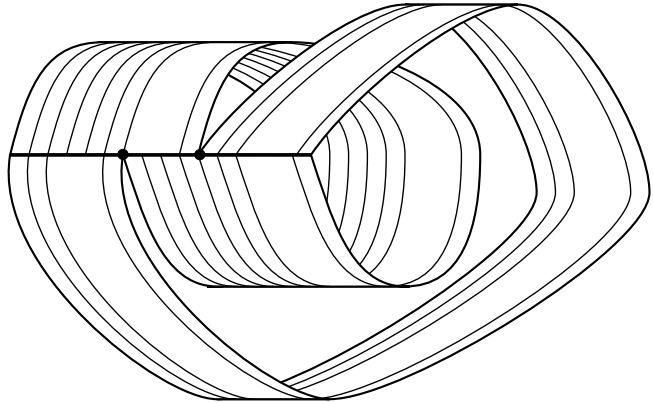}
\end{center} 
\caption{\label{fig:iet} Vertical foliation of the mapping torus of an
  interval exchange transformation.  This foliation has two
  singularities (in bold), each of Poincar\'e-Lefschetz index $-\frac
  12$ and local $\CQ$-index $1$. The $\CQ$-index of the foliation is
  $2$. The surface is a torus with one boundary component: its Euler
  Characteristic is $-1$.}
\end{figure}

Interval exchange transformations provide examples of such
foliated surfaces. Indeed, the mapping torus of an interval
exchange transformation is a surface (with boundary), naturally
foliated by the vertical direction, as in Figure~\ref{fig:iet}.
The transverse measure of the foliation is given by the Lebesgue
measure of the interval. We define, in this case, the index of
the interval exchange transformation as the index defined above
for this foliation and its dual $\R$-tree.

\subsection{Geometric trees}

This surface situation has been generalized (see for instance
\cite{best-survey}). Let us consider a finite family of partial
isometries of an interval (or a finite number of intervals, or even a
finite tree). The suspension of these partial isometries gives rise to
a 2-complex $B$ (which is not a surface in general), naturally
foliated by the vertical direction. As previously, the foliation can
be lifted to the universal cover of $B$, and the space of leaves of
this foliation is an $\R$-tree with an action of $\pi_1(B)$ by
isometries. A tree obtained by such a construction is called a
geometric tree.  In this situation, we can define two local indices:
one for the tree, using the valence of branch points, and one for the
foliation, using the number of ends of singular leaves. We would like
to stress that, contrary to the case of a foliation on a surface,
these two indices do not agree, not even locally.

The first index has been introduced by Gaboriau and
Levitt~\cite{gl-rank}. In this paper we call it the geometric index,
and denote it by $\igeo(T)$. It is defined using the valence of the
branch points, of the $\R$-tree $T$, with an action of the free group
by isometries:
\[ 
\igeo(T)=\sum_{[P]\in T/\FN}\igeo(P). 
\]
where the local index of a point $P$ in $T$ is
\[
\igeo(P)=\#(\pi_0(T\smallsetminus\{P\})/\Stab(P))+2\,\rank(\Stab(P))-2.
\]
Gaboriau and Levitt \cite{gl-rank} proved that the geometric index of
a geometric tree is equal to $2N-2$ and that for any tree in the
compactification of Outer Space $\barCVN$ the geometric index is
bounded above by $2N-2$. Moreover, they proved that the trees in
$\barCVN$ with geometric index equal to $2N-2$ are precisely
the geometric trees.

The second index is defined from the number of ends of singular
leaves. To our knowledge it has never been studied in its own right,
although Gaboriau~\cite{gab-bouts} gives a lot of relevant
insights. In particular Gaboriau \cite[Theorem VI.1]{gab-bouts} gives
partial results to bound this index.

\subsection{The $\CQ$-index of an $\R$-tree}

Let $T$ be an $\R$-tree in the boundary of Outer Space with dense
orbits. We denote by $\bar T$ its metric completion, $\partial T$ its
Gromov-boundary at infinity.  The set $\hat T=\bar T\cup\partial T$
equipped with the observers' topology (a slight weakening of the
metric topology, see \cite{chl2}) is a compact set denoted $\Tobs$.

Let $P$ be a point in $T$.  The map $\CQ: \partial\FN \rightarrow
\Tobs$ is the unique continuous extension (see \cite{chl2}) of the
map
\[
\begin{array}{rcl}
\FN & \rightarrow & T\\
u & \mapsto & u\cdot P.
\end{array}
\]
It does not depend on the choice of the point $P$.

The map $\CQ$ can be easily understood in the special case of a tree $T$ 
dual to a foliation on a 2-complex $B$.  Each leaf of the
foliation of $\tilde{B}$ is a point of the dual tree $T$.  A half-leaf
of the foliation of $\tilde{B}$ defines a point
$X\in\partial\FN=\partial\pi_1(B)$ and $\CQ(X)$ is the point of $T$
defined by the leaf.

The general definition of the $\CQ$-index of the tree $T$ is given as
follows:
\[ 
\iq(T)=\sum_{[P]\in \hat T/\FN}\max(0;\iq(P)). 
\] 
where the local index of a point $P$ in $T$ is:
\[ 
\iq(P)=\#(\CQ\inv(P)/\Stab(P))+2\,\rank(\Stab(P))-2. 
\]

Levitt and Lustig~\cite{ll-north-south} proved that points in $\partial T$ have
exactly one pre-image by $\CQ$ (see Proposition~\ref{prop:Qonto}).
Thus, only points in $\bar T$ contribute to the $\CQ$-index of $T$.

Our main result states that the $\CQ$-index of an $\R$-tree in the
boundary of Outer Space is bounded above by $2N-2$. This answers the
question of Levitt and Lustig \cite[Remark~3.6]{ll-north-south}
whether the map $\CQ:\partial\FN\to\hat T$ has finite fibers (in the
case where the action is free).

In \cite{chl1-II} the dual lamination of $T$ is defined: it is the set
of pairs $(X,Y)$ of distinct points in the boundary $\partial\FN$ such
that $\CQ(X)=\CQ(Y)$. More conceptually, the $\CQ$-index of $T$ can
indeed be understood in a more general context as that of its dual
lamination.

The limit set $\Omega$ is the subset of $\bar T$ which consists of
points with at least two pre-images by the map $\CQ$.

We also describe the trees such that $\iq(T)=2N-2$: these are the
trees such that all points of $T$ have at least two pre-images by
$\CQ$.

\begin{thmsurface}
  Let $T$ be an $\R$-tree in the boundary of Outer Space with dense
  orbits. The $\CQ$-index is maximal: $\iq(T)=2N-2$ if and only if $T$
  is contained in the limit set $\Omega$.
\end{thmsurface}

An $\R$-tree dual to a foliation on a surface with boundary of
negative Euler characteristic, has maximal $\CQ$-index. We call trees
with maximal $\CQ$-index trees of surface type.

\subsection{Compact systems of isometries}

A traditional strategy to study a tree in the boundary of Outer
Space is 
\begin{enumerate}
\item describe any geometric tree by a system of isometries on a finite
tree (or even a multi-interval) and then use the Rips Machine;  
\item approximate any tree by a sequence of geometric trees.
\end{enumerate}
In particular Gaboriau and Levitt~\cite{gl-rank} proved in this way
that the geometric index of any tree in Outer Space is bounded
above by $2N-2$.

In \cite{chl4} a new approach was proposed: to describe an $\R$-tree
by a system of isometries on a compact $\R$-tree (rather than on a
finite tree). The point here is that any tree $T$ in the
compactification of Outer Space can be described by a system of
isometries on a compact $\R$-tree: $S_A=(K_A,A)$ (where $A$ is a basis
of $\FN$ and $K_A$ is a compact subtree of $\bar T$). This system of
isometries encodes all of the original tree $T$ (together with the
action of $\FN$). An index is defined in Section~\ref{sec:indexsi} for
any such system of isometries.

\begin{thmsiiq} 
  Let $T$ be an $\R$-tree with a very small, minimal action of $\FN$
  by isometries with dense orbits. The $\CQ$-index of $T$ and the
  index of the induced system of isometries $\CK_\CA=(K_\CA,\CA)$, for
  any basis $\CA$, are equal:
\[
\iq(T)=i(\CK_\CA). 
\] 
\end{thmsiiq}

The computation of the index of a tree is thus achieved  by  computing
the index of a system of isometries. We study system of isometries by
themselves in Sections~\ref{sec:si}, \ref{sec:ripsmachine} and
\ref{sec:compindexsi}.

We improve the classical Rips Machine (see \cite{glp-rips,bf-stable})
to work in the context of systems of isometries on compact $\R$-trees
(or forests). The Rips Machine applied to a system of isometries
returns a new system of isometries obtained by erasing parts of the
supporting forest. To each system of isometries we associate a finite
graph $\Gamma$, the index of which is given by the Euler
characteristic.  We study the effect of the Rips Machine on this
graph: the Rips Machine decreases the index of the graph $\Gamma$.

Iterating the Rips Machine infinitely many times, the sequence of
associated graphs $\Gamma$ has a limit $\hat\Gamma$. The index of
$\hat\Gamma$ is bounded above by the decreasing sequence of
indices. We prove that the index of the limit graph $\hat\Gamma$ is
equal to the index of the system of isometries. In fact, in the case
of a Levitt system of isometries the graph $\hat\Gamma$ can be viewed
as a geometric realization of the dual lamination of the system of
isometries. We obtain

\begin{thmmainsi} 
The index of a system of isometries $\CK=(\CF,\CA)$ with independent
generators is finite and bounded above by the index of the associated
graph $\Gamma$.
\end{thmmainsi}

The above Theorem~\ref{thm:surface} follows from our characterization
of systems of isometries with maximal index:

\begin{thmmaxindex}
Let $\CK=(\CF,\CA)$ be a reduced system
of isometries with independent generators, let $\Gamma$ be its
associated graph, and $\hat\Gamma$ be its limit graph. The following are equivalent
\begin{enumerate} 
\item The system of
isometries $\CK$ has maximal index,
\item The graph $\hat\Gamma$ is
finite, 
\item The Rips Machine, starting
from $\CK$, halts after finitely many steps.
\end{enumerate}
\end{thmmaxindex}

\bigskip

\noindent\textbf{Acknowledgment:} {\em We thank Martin Lustig for his
  constant interest in our work.  

We are grateful to Vincent Guirardel and Gilbert Levitt for
introducing us to mixing properties of trees.
}

\section{Systems of isometries}\label{sec:si}

\subsection{Definitions}\label{subsec:sidef}

We collect in this Section basic facts from \cite{chl4}.

An \textbf{$\R$-tree}, $(T,d)$ is a metric space such that for
any two points $P,Q$ in $T$, there exists a unique arc $[P;Q]$
between them and this arc is isometric to the segment
$[0;d(P,Q)]$.

A \textbf{compact forest} $\CF$ is a metric space with finitely
many connected components each of which is a compact $\R$-tree.

A \textbf{partial isometry} of a compact forest $\CF$ is an
isometry $a:K\to K'$, between two compact subtrees $K$ and $K'$
of $\CF$. The \textbf{domain} of $a$ is $K$, its \textbf{range}
is $K'$. The partial isometry $a$ is \textbf{non-empty} if its
domain is non-empty. The domain (and the range) of a partial
isometry needs not be a whole connected component of $F$. A
\textbf{system of isometries} $\CK=(\CF,\CA)$ consists of a
compact forest $\CF$ and of a finite set $\CA$ of non-empty
partial isometries of $\CF$.

To such a system of isometries $\CK$ we associate the oriented
graph $\Gamma$ which has the connected components of $\CF$ as
vertices and the non-empty partial isometries of $\CA$ as
oriented edges. The edge $a\in\CA$ starts at the connected
component of $\CF$ which contains its domain, and ends at the
connected component of $\CF$ which contains its range.

We regard the reverse edge $a\inv$ of the edge $a\in\CA$ as the
inverse partial isometry $a\inv$ of $a$. A reduced path $w$ in
the graph $\Gamma$, given as a sequence of edges $w=z_1\cdots
z_n$ with $z_i\in\CA^{\pm 1}$ (such that $z_{i+1}\neq
{z_i}\inv$), defines a (possibly empty) partial isometry, also
denoted by $w$: the composition of partial
isometries $z_1 \circ z_{2} \circ \cdots \circ z_n$. We write
this \textbf{pseudo-action} on $\CF$ on the right, i.e. \[ P (u
\circ v) = (Pu)v \] for all point $P \in \CF$ and for all path
$uv$ in $\Gamma$.

The \textbf{pseudo-orbit} of a point $P$ in $\CF$ is the subset
of $\CF$ which can be reached from $P$: 
\[ 
\{P.w\ |\ w\mbox{ reduced path in }\Gamma, P\in\dom(w)\}. 
\]

A reduced path $w$ in $\Gamma$ is \textbf{admissible} if it is
non-empty as a partial isometry of $\CF$.

An infinite reduced path $X$ in $\Gamma$ is \textbf{admissible}
if all its subpaths are admissible. The domains of the initial
subpaths of $X$ are nested compact subtrees, their intersection
is the \textbf{domain} of $X$, denoted by $\dom(X)$.

A bi-infinite reduced path $Z$ in $\Gamma$ is
\textbf{admissible} if all its subpaths are admissible. A
bi-infinite reduced path, $Z=\cdots
z_{-2}z_{-1}z_{0}z_{1}z_{2}\cdots$, has two halves which are
infinite reduced paths: 
\[ 
Z^+=z_1z_2\cdots,\quad Z^-={z_{0}}\inv{z_{-1}}\inv{z_{-2}}\inv\cdots.
\] 
The \textbf{domain} of $Z$ is the intersection of the domains of its
two halves. Equivalently, a bi-infinite reduced path $Z$ in $\Gamma$
is admissible if and only if its domain is non-empty.

The \textbf{limit set} $\Omega$ of a system of isometries
$\CK=(\CF,\CA)$ is the set of elements of $K$ which are in the
domain of a bi-infinite admissible reduced path in $\Gamma$. The
limit set is the place where the dynamics of the system of
isometries concentrates. Alternatively, $\Omega$ is the largest
subset of $\CF$ such that for each $P\in\Omega$ there exists at
least two partial isometries $a,b\in\CA^{\pm 1}$ with $P.a$ and
$P.b$ in $\Omega$.

A system of isometries $\CK=(\CF,\CA)$ has \textbf{independent
generators} (compare Gaboriau \cite{gab-indgen} and \cite{chl4})
if the domain of any admissible infinite reduced path $X$ in
$\Gamma$ consists of exactly one point which we denote by
$\CQ(X)$: 
\[ 
\dom (X)=\{ \CQ(X)\}. 
\] 
In this case, the domain of a bi-infinite admissible reduced path $Z$
in $\Gamma$ also consists of exactly one point which we also denote by
$\CQ(Z)$.

If $\CK$ has independent generators we have: 
\[ 
\Omega=\{ P\in K | \exists Z\mbox{ bi-infinite admissible, }
\CQ(Z)=P\}.
\]

The \textbf{restriction} of a partial isometry $a$ to the
compact $\R$-tree $K$ (at the source) is the (possibly empty)
partial isometry ${}_{K\rceil}a$ which is defined for each $P\in
K\cap\dom (a)$. The restriction of $a$ to the compact $\R$-tree
$K'$ (at the destination) is the (possibly empty) partial
isometry $a_{\lceil K'}$ which is defined for each $P\in\dom
(a)$ such that $Pa\in K'$. The restriction of the partial
isometry to the compact $\R$-trees $K$ and $K'$ is the (possibly
empty) partial isometry ${}_{K\rceil}a_{\lceil K'}$ which is
defined for each $P\in K\cap\dom (a)$ such that $P.a$ is in
$K'$.

\subsection{Index of a Graph}

We denote by $V(\Gamma)$ the set of vertices of a graph $\Gamma$
and by $E(\Gamma)$ its set of edges.

For a vertex $x$ of a graph $\Gamma$ the \textbf{valence}
$v_\Gamma(x)$ of $x$ is the number of edges incident to $x$. The
\textbf{index} $i_\Gamma(x)=v_\Gamma(x)-2$ of $x$ is its valence
minus two.

The \textbf{index} $i(\Gamma)$ of a finite connected graph
$\Gamma$ is: 
\[ 
\begin{array}{rcl}i(\Gamma)&=&\max(0;\sum_{x\in V(\Gamma)}i_\Gamma(x))\\ 
&=&\max(0;2(\#E(\Gamma)-\#V(\Gamma)))\\
&=&\max(0;-2\chi(\Gamma))\\ &=&\max(0;2\,\rank(\pi_1(\Gamma))-2)
\end{array} 
\] 
where $\chi(\Gamma)$ is the Euler characteristic
of $\Gamma$. The index $i(\Gamma)$ is a homotopy invariant of
the graph $\Gamma$.

\begin{figure}[h]
\centering
\input{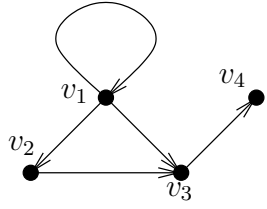}
\caption{\label{fig:index}A graph of index 2, with vertices of indices
  $i(v_1)=2$, $i(v_2)=0$, $i(v_3)=1$ and $i(v_4)=-1$.}
\end{figure}

The \textbf{index} of a finite graph $\Gamma$ is the sum of the
indices of its connected components.

The \textbf{core} of a graph $\Gamma$ is the largest subgraph of
$\Gamma$ without vertices of valence $0$ and $1$ (see \cite{gers}), we
denote it by $\core(\Gamma)$. The core of a graph is the union of all
bi-infinite reduced paths. Alternatively, if $\Gamma$ is finite, the
subgraph $\core(\Gamma)$ is obtained by recursively erasing vertices
of valence $0$ and the vertices of valence $1$ together with their
incident edges. The core of a graph may be empty: this is the case,
for instance, when the graph is a union finite of trees.

The index of a finite graph can be computed using its core with
the easier formula: 
\[ 
i(\Gamma)=i(\core(\Gamma))=\sum_{x\in
V(\core(\Gamma))}i_{\core(\Gamma)}(x).
\]

We use this formula to define the index of an infinite graph
$\Gamma$: The core  does not have vertices of
valence $0$ and $1$ and thus we can compute the above (possibly
infinite) non-negative sum.

For a connected (possibly infinite) graph $\Gamma$ the index is
thus 
\[ 
\begin{array}{rcl}i(\Gamma)&=&i(\core(\Gamma))\\
&=&\max(0;\#\partial\Gamma+2\,\rank(\pi_1(\Gamma))-2)
\end{array} 
\] 
where $\partial\Gamma$ is the set of ends of $\Gamma$. If $\Gamma$ is
not connected we sum the above value for each of its connected
components.

A \textbf{morphism of graphs} $\tau:\Gamma'\to\Gamma$ maps
vertices to vertices, edges to edges and respects incidence.

We will need the following Lemma in our proofs:

\begin{lem}\label{lem:indexgraph} 
Let $\tau:\Gamma'\to\Gamma$ be
a morphism between two finite graphs. Assume that $\tau$ is
injective on edges. Then the index of $\Gamma'$ is smaller or
equal to the index of $\Gamma$: 
\[ 
i(\Gamma')\leq i(\Gamma). 
\]
\end{lem} 
\begin{proof} 
For each vertex $x$ of $\Gamma$, the set
of edges incident to the vertices in the fiber $\tau\inv(x)$
injects in the set of edges incident to $x$. Thus 
\[
\sum_{x'\in\tau\inv(x)}i_{\Gamma'}(x')\leq i_{\Gamma}(x). 
\]
Moreover, $\tau$ maps the core  of $\Gamma'$
inside the core  of $\Gamma$.

We get 
\[ 
i(\Gamma')=i(\core(\Gamma'))=\sum_{x\in
V(\core(\Gamma'))}i_{\core(\Gamma')}(x)\leq \sum_{x\in
V(\core(\Gamma))}i_{\core(\Gamma)}(x)=i(\core(\Gamma))=i(\Gamma).
\] 
\end{proof}

In particular if $\Gamma'$ is a subgraph of a finite graph
$\Gamma$ 
\[ 
i(\Gamma')\leq i(\Gamma). 
\]

\subsection{Index of a system of isometries}\label{sec:indexsi}

Let $\CF$ be a compact forest and $\CK=(\CF,\CA)$ be a system of
isometries. Let $\Omega$ be the limit set of $\CK$.

For a point $P$ in $\CF$, we define its \textbf{index} by 
\[
i_\CK(P)=\#\{a\in\CA^{\pm 1}\ |\ P.a\in\Omega\}-2. 
\] 
By definition of the limit set, for any point $P$ in $\Omega$, there
exists a bi-infinite reduced admissible path $Z=\cdots
z_{-1}z_0z_1\cdots$ in $\Gamma$ such that $P\in\dom(Z)$. The edges
$z_1$ and ${z_0}\inv$ send $P$ inside $\Omega$, and thus the index of
$P$ is greater or equal to $0$: $i_\CK(P)\geq 0$.

Conversely, if the index of a point $P$ in $F$ is non-negative:
$i_S(P)\geq 0$, then there exists two elements $a,b\in\CA^{\pm
1}$, such that $P.a$ and $P.b$ are in the limit set $\Omega$. As
$P.a$ is in $\Omega$ there exists a bi-infinite reduced
admissible path $Z=\cdots z_{-1}z_0z_1\cdots$ in $\Gamma$ such
that $P.a\in\dom(Z)$. Up to replacing $Z$ by $Z\inv$ (the same
bi-infinite path with reversed orientation), we assume that
$z_1\neq a\inv$. Symmetrically there exists a bi-infinite
reduced admissible path $Z'=\cdots z'_{-1}z'_0z'_1\cdots$ in
$\Gamma$ such that $P.b\in\dom(Z')$ and $z'_0\neq b\inv$. We get
that $Z''=\cdots z'_{-1}z'_0b\cdot a z_1z_2\cdots$ is a
bi-infinite reduced path in $\Gamma$ with $P\in\dom(Z'')$ and
thus that $P$ is in the limit set $\Omega$: 
\[ 
P\in\Omega\iff i_S(P)\geq 0. 
\]

The \textbf{index} of $\CK$ is defined by 
\[
i(\CK)=\sum_{P\in\CF}\max(0; i_\CK(P))=\sum_{P\in\Omega}i_\CK(P).
\]
As there is a $\max$ in the first sum, and by the above equivalence,
both sums are non-negative and are well defined possibly as $+\infty$.

The main result of this paper can now be stated

\begin{thm}\label{thm:mainsi} 
The index of a system of isometries $\CK=(\CF,\CA)$ with independent
generators is finite and bounded above by the index of the associated
graph $\Gamma$.
\end{thm}

\subsection{Cayley graphs}\label{sec:cayley}

Let $\CK=(\CF,\CA)$ be a system of isometries and let $\Gamma$
be its associated graph. Let $P$ be a point in $\CF$ and $K$ be
the connected component of $\CF$ which contains $K$. Let
$\Gamma_0$ be the connected component of $\Gamma$ which contains
$P$ and $\tilde\Gamma_0$ its universal cover.

The \textbf{trajectory tree} of $P$ is the smallest subtree
$\Tt(P)=\Tt(P,S)$ of $\tilde\Gamma_0$ which contains all the
admissible paths $w$ based at $K$ such that $P$ is in the domain
of $w$.

Let $\Stab(P)$ be the subgroup of the fundamental group
$\pi_1(\Gamma,K)$ of the graph $\Gamma$ based at $K$ of
admissible paths $w$ such that $P.w=P$. The group $\Stab(P)$ is
a free group that acts on the tree $\Tt(P)$.

The \textbf{Cayley graph} $\Cayl(P)=\Cayl(P,S)$ of $P$ is the
quotient of $\Tt(P)$ by the action of the stabilizer $\Stab(P)$
(compare Gaboriau \cite{gab-bouts}). Vertices of the Cayley
graph of $P$ are in one-to-one correspondence with the elements
of the pseudo-orbit of $P$ in $\CF$ under the pseudo-group of
isometries. The vertices of the core of the
Cayley graph of $P$ are in one-to-one correspondence with the
intersection $\omega(P)$ of the pseudo-orbit of $P$ and the
limit set $\Omega$.

The index of a vertex $P'$ in $\core(\Cayl(P))$ is
equal to the index $i_S(P')$ of the point $P'$ for the system of
isometries $S$. Thus, we get that the index of the core of the
Cayley graph $i(\Cayl(P))$ is equal to the contribution of the
pseudo-orbit of $P$ to the index of $\CK$:
\[
i(\Cayl(P))=i(\core(\Cayl(P))=\sum_{P'\in\omega(P)}i_\CK(P').
\]
Adding up, for all pseudo-orbits $[P]$ we get
\[
i(S)=\sum_{[P]} i(\Cayl(P)).
\]

\section{Rips Machine}\label{sec:ripsmachine}

\subsection{Elementary step}

Let $\CK=(\CF,\CA)$ be a system of isometries on a compact
forest $\CF$. One step of the Rips Machine produces a new system
of isometries $\CK'=(\CF',\CA')$ defined as follows.

The forest $\CF'$ is the set of all elements of $\CF$ which are
in the domains of at least two distinct partial isometries in
$\CA^{\pm 1}$:
\[
\CF'=\{P\in\CF\ |\ \exists a\neq b\in \CA^{\pm 1}, P\in \dom (a)\cap
\dom (b)\}.
\]
The set $\CF'$ has finitely
many connected components which are compact $\R$-trees because
it is the finite union of all possible intersections $\dom
(a)\cap \dom (b)$ for all pairs of distinct elements $a,b$ of
$\CA^{\pm 1}$.

For each partial isometry $a\in\CA$ and for each pair of
connected components $K_0,K_1$ of $\CF'$, we consider the
partial isometry, $a'={}_{K_0\rceil}a_{\lceil K_1}$, which is
the restriction of $a$ to $K_0$ and $K_1$. The finite set $\CA'$
consists of all such non-empty partial isometries
$_{K_0\rceil}a_{\lceil K_1}$ of $\CF'$.

An elementary step of the Rips Machine gives rise to a map
$\tau$ from the graph $\Gamma'$ associated to the resulting
system of isometries $\CK'$, to the original graph $\Gamma$. A
vertex $K'$ of $\Gamma'$ is a connected component of $\CF'$ and
it is mapped by $\tau$ to the connected component $\tau(K')$ of
$\CF$ which contains $K'$. Similarly an edge $a'$ of $\Gamma'$
is a non-empty partial isometry $a'={}_{K_0\rceil}a_{\lceil
K_1}$ and it is mapped by $\tau$ to the original partial
isometry $a$. The map $\tau$ is a morphism of oriented graphs.

If $w$ is an admissible reduced path in $\Gamma'$, the domain of
$\tau(w)$ contains the domain of $w$ and $\tau(w)$ is an
admissible path in $\Gamma$. Moreover, the image $\tau(w)$ of an
admissible reduced path $w$ in $\Gamma'$ is a reduced path of
$\Gamma$. Finally, if a bi-infinite reduced path $Z$ in $\Gamma$
is admissible, then its domain is contained in $\CF'$, which
leads to the following

\begin{prop}\label{prop:limitsetstepRips} 
Let $\CK$ be a system of isometries and $\CK'$ be the result of the
Rips Machine. Then the limit sets and the indices of $\CK$ and $\CK'$
are equal:
\[
\Omega=\Omega'\mbox{ and }i(S)=i(S'). \qed
\]
\end{prop}

\subsection{Indices through the Rips Machine}

As explained previously, the Rips Machine defines a new system
of isometries $\CK'=(\CF',\CA')$ starting from a system of
isometries $\CK=(\CF,\CA)$ by erasing the subset $E$ of the
forest $F$ which consists of points which belongs to at most one
domain of partial isometries of $\CA^{\pm 1}$.

To better understand the Rips Machine we decompose its
elementary step into a finite sequence of elementary moves.
Instead of erasing $E$ in one step we successively erase subsets
$E_i$ of $E$. This gives us a finite sequence of system of
isometries starting from $\CK$ and ending at $\CK'$. The
successive systems of isometries of this sequence differ by an
elementary move.

As $\CF$ is a compact forest, the set $E$ may have infinitely
many connected components. We first describe a preliminary move
which erases all the ``peripheral'' ones. Then we are left with
a finite forest to erase, which we erase in finitely many
elementary moves.

This decomposition of the Rips Machine is used in the next
Proposition to prove that the index of the associated graphs is
decreasing.

\begin{prop}\label{prop:indexripsmachine} 
Let $\CK=(\CF,\CA)$ be a system of isometries and $\CK'=(\CF',\CA')$
be the output of the Rips Machine. Let $\Gamma$ and $\Gamma'$ be the
associated graphs.

Then the index $i(\Gamma')$ is smaller or equal than the index
$i(\Gamma)$.
\end{prop}
\begin{proof} Let $E$ be the part of the
forest $\CF$ which is erased by the Rips Machine:
\[
E=\CF\smallsetminus\CF'=\{P\in\CF\ |\ \#\{a\in\CA^{\pm 1}\ |\
P\in\dom (a)\}\leq 1\}.
\]
Let $E_C$ be the subset of $E$ which is in the convex hull of $\CF'$:
\[
E_C=\{P\in E\ |\ \exists
Q,R\in\CF', P\in [Q,R]\},
\]
and let $E_0=E\smallsetminus E_C$ be the complement of $E_C$ in $E$.

Let $\CF_0=\CF\smallsetminus E_0=\CF'\cup E_C$: $\CF_0$ is the
convex hull of the connected components of $\CF'$ in $\CF$:
\[
\CF_0=\{P\in \CF\ |\ \exists Q,R\in\CF', P\in [Q,R]\}.
\]
Thus $\CF_0$ has finitely many connected components, each of which is
a compact $\R$-tree: $\CF_0$ is a compact forest.

Let $\CA_0$ be the set of all non-empty restrictions of partial
isometries of $\CA$ to partial isometries of $\CF_0$. Let
$\Gamma_0$ be the graph associated to the system of isometries
$\CK_0=(\CF_0,\CA_0)$ and let $\tau_0:\Gamma_0\to\Gamma$ be the
graph morphism defined as before. As $E_0$ does not split
connected components of $\CF$, the map $\tau_0$ is injective and
thus by Lemma~\ref{lem:indexgraph},
\[
i(\Gamma_0)\leq i(\Gamma)
\]
(we note that this inequality can be strict if $\tau_0$ is
not onto).

As $\CF'$ has finitely many connected components, $E_C$ is a
finite union of finite open arcs of the form $]P;Q[$ where $P$
and $Q$ are two points of $\CF'$.

\begin{lem}\label{lem:arcs} 
The erased part $E_C$ can be
decomposed to get a partition
\[
E_C=\alpha_1\uplus \alpha_2
\uplus\cdots\uplus \alpha_n
\]
where each $\alpha_i$ is an open
arc such that the number of connected components of $\CF_i$ is
exactly one plus the number connected components of $\CF_{i-1}$,
where for each $i=1,\ldots,n$, we let
$\CF_i=\CF_{i-1}\smallsetminus \alpha_i$. 
\end{lem}
\begin{proof} 
We recursively define $\alpha_i$ by choosing a
connected component $C$ of $\CF_{i-1}$ which contains at least
two connected components of $\CF'$. Then we choose a connected
component $K$ of $\CF'$ contained in $C$ and which is not
contained in the convex hull $C'$ of $(C\cap\CF')\smallsetminus
K$ in $C$. We choose the open arc $\alpha_i$ that joins $K$ and
$C'$. Removing $\alpha_i$ from $\CF_{i-1}$ splits the connected
component $C$ of $\CF_{i-1}$ into two new connected components:
$K$ and $C'$. 
\end{proof}

Let $\CA_i$ be the set of non-empty restrictions of partial
isometries in $\CA$ to $\CF_i$ and let $\Gamma_i$ be the graph
associated to the system of isometries $\CK_i=(\CF_i,\CA_i)$. As
before we get graph morphisms
$\tau_i:\Gamma_{i}\to\Gamma_{i-1}$. We observe that the last
system of isometries is the output of the Rips Machine:
$\CK_n=\CK'$. The map $\tau$ factors through the graphs
$\Gamma_i$: $\tau=\tau_0\circ\tau_1\circ\cdots\circ\tau_n$.

We now proceed to prove that for each $i=1,\ldots,n$, the index
of $\Gamma_i$ is lower or equal to the index of $\Gamma_{i-1}$.
This will conclude the proof.

Going from $\CK_{i-1}$ to $\CK_{i}$ corresponds to one of the
following \textbf{elementary moves} of the Rips Machine. For
each $i$ removing the arc $\alpha_i$ from $\CF_{i-1}$ has one of
the following two effects on the graph $\Gamma_{i-1}$:
\begin{enumerate} 
\item \underline{Split a vertex}: $\alpha_i$
joins two connected components, $K$ and $K'$, of $\CF_i$, and no
partial isometry in $\CA_{i-1}^{\pm 1}$ is defined
simultaneously on $K$ and $K'$. Then the map $\tau_i$ is
injective on edges and, applying Lemma~\ref{lem:indexgraph}, the
index of $\Gamma_i$ is smaller or equal to the index of
$\Gamma_{i-1}$. 

\begin{figure}[h]
\centering
\input{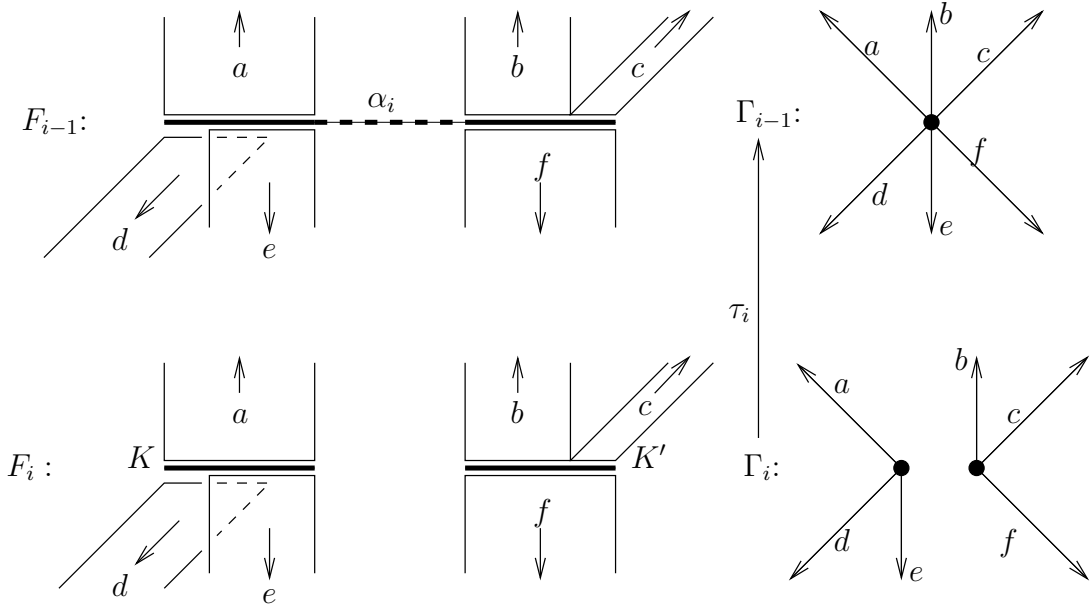}
\caption{\label{fig:split-vertex}Split a vertex move}
\end{figure}

\item \underline{Split an edge}: $\alpha_i$
joins two connected components $K$ and $K'$ of $\CF_i$, and a
partial isometry $a$ in $\CA_{i-1}^{\pm 1}$ is defined on $K$
and $K'$ (and thus its domain $\dom(a)$ contains $\alpha_i$). By
definition of $E$ no other partial isometry is defined on
$\alpha_i$, in particular the range of $a$ is contained in
$\CF_i$. The graph morphism $\tau_i$ maps the two distinct
vertices $K$ and $K'$ to the same vertex, $K\cup\alpha_i\cup
K'$, of $\Gamma_{i-1}$, it maps the two edges $a'={}_{K\rceil}a$
and $a''={}_{K'\rceil}a$ which are the restrictions of $a$ to
$K$ and $K'$ to the same edge $a$ of $\Gamma_{i-1}$. On all
other vertices and edges, $\tau_i$ is one-to-one. Thus, $\tau_i$
is a homotopy equivalence and the index of $\Gamma_i$ is equal
to the index of $\Gamma_{i-1}$. Indeed, $\tau_i$ is a folding in
Stallings' terminology \cite{stall-graphs}. 

\begin{figure}[h]
\centering
\input{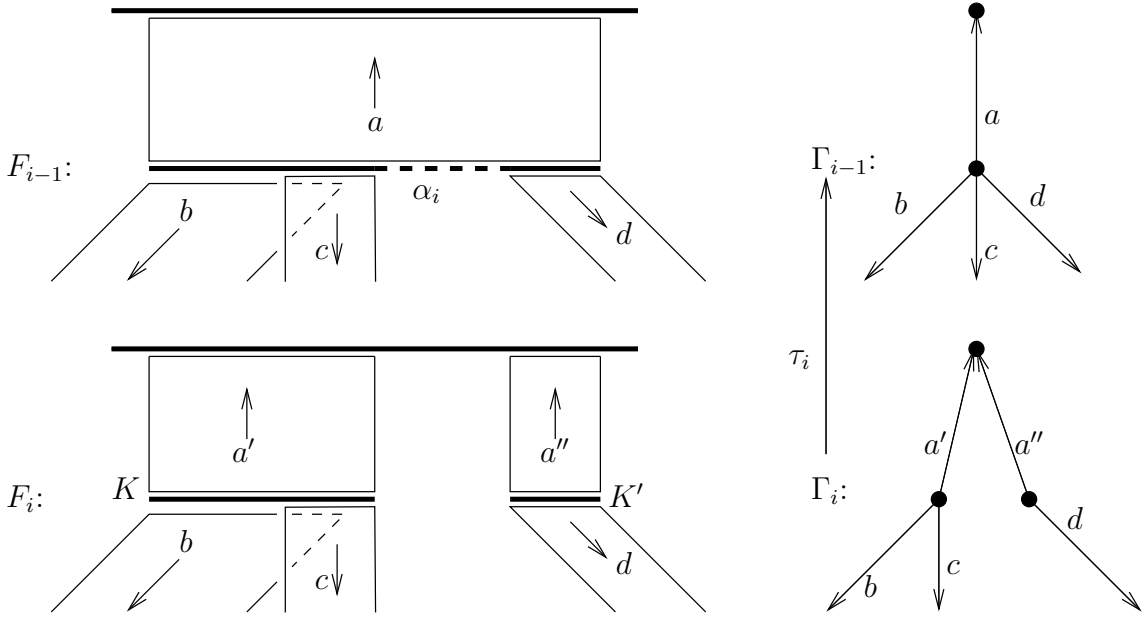}
\caption{\label{fig:split-edge}Split an edge move}
\end{figure}

\end{enumerate}
\end{proof}

\subsection{Iterating the Rips Machine}

Let $\CK_0=(\CF_0,\CA_0)$ be a system of isometries. By
repeatedly applying the Rips Machine we get a sequence
$\CK_n=(\CF_n,\CA_n)$ of systems of isometries. Of course the
Rips Machine may halt after some time, that is to
say, we do not exclude that $\CK_{n+1}=\CK_n$ for $n$ big
enough.

We also get graph morphisms $\tau_n$ from the graph
$\Gamma_{n+1}$ associated to $\CK_{n+1}$ to the graph $\Gamma_n$
associated to $\CK_n$. Indeed, a partial isometries
$a_n\in\CA_n$ is the restriction of the partial isometry
$a=\tau_0\circ\tau_1\cdots\tau_{n-1}(a_n)$ in $\CA_0$ to
connected components of the compact forest $\CF_n$.

\begin{lem}\label{lem:limitrips} 
The intersection $\Omega$ of
the nested sequence $(\CF_n)_{n\in\N}$ of compact subsets of
$\CF_0$ is equal to the limit set $\Omega_0$ of the system of
isometries $\CK_0$. 
\end{lem} 
\begin{proof} 
By Proposition~\ref{prop:limitsetstepRips}, at each step $n$ the limit
sets $\Omega_{n+1}$ and $\Omega_n$ of the corresponding system of
isometries are equal. In particular $\Omega_0$ is contained in $\CF_n$
at each step $n$ and thus in the nested intersection.

Conversely, let $P$ be a point in the nested intersection
$\Omega=\cap_{n\in\N}\CF_n$. For any $n\in\N$, $P$ belongs to
$\CF_{n+1}$ and by definition of the Rips Machine, there exists
at least two distinct partial isometries $a_{n}$ and $b_{n}$
defined at $P$ in ${\CA_n}^{\pm 1}$. Thus $P.a_n$ and $P.b_n$
are in $\CF_n$. Up to passing to a subsequence we can assume
that $a_n$ and $b_n$ are the restrictions of two fixed and
distinct partial isometries $a$ and $b$ in $\CA_0^{\pm 1}$. We
get that for all $n$, $P.a$ and $P.b$ are in $\CF_n$, which
proves that $P.a$ and $P.b$ are also in the nested intersection
$\Omega=\cap_{n\in\N}\CF_n$. The set $\Omega$ is a subset of
$\CF_0$ such that for any point $P$ in $\Omega$, there exists at
least two distinct partial isometries $a$ and $b$ in $\CA_0^{\pm
1}$ such that $P.a$ and $P.b$ are also in $\Omega$. This proves
that $\Omega$ is contained in the limit set $\Omega_0$ of the
initial system of isometries $\CK_0$. 
\end{proof}

The \textbf{limit graph} $\hat\Gamma$ of the system of
isometries $\CK_0=(\CF_0,\CA_0)$ is the (possibly infinite)
graph whose vertices are the connected components of the limit
set $\Omega$ and whose edges are all possible restrictions of
partial isometries in $\CA_0$ to connected components of
$\Omega$. We denote by $\hat\tau_n:\hat\Gamma\to\Gamma_n$ the graph
morphism that maps a connected component $C$ of $\Omega$ to the
connected component of $\CF_n$ that contains $C$, and which maps
an edge $e$ of $\hat\Gamma$ to the partial isometry
$a_n\in\CA_n$ of which it is a restriction.

From the previous Lemma and from the definition of an inverse
limit, we deduce:

\begin{lem}\label{lem:inverselimit} 
The limit graph $\hat\Gamma$
of a system of isometries, $\CK_0=(\CF_0,\CA_0)$, is the inverse
limit of the sequence of graphs $(\Gamma_n)_{n\in\N}$ (together
with the sequence of maps $(\tau_n)_{n\in\N}$) associated to the
sequence of systems of isometries $(\CK_n)_{n\in\N}$ obtained
from $\CK_0$ by iterating the Rips Machine. \qed
\end{lem}

By definition of $\Omega$, the graph $\hat\Gamma$ does not have
vertices of valence $0$ or $1$, and we defined its
\textbf{index} as the non-negative sum
\[
i(\hat\Gamma)=\sum_{x\in V(\hat\Gamma)}i_{\hat\Gamma}(x).
\]
Recall that, by Proposition~\ref{prop:indexripsmachine}, the
sequence of indices $(i(\Gamma_n))_{n\in\N}$ is decreasing.

\begin{prop}\label{prop:hatgamma} Let $\hat\Gamma$ be the limit
graph of a system of isometries $\CK_0=(\CF_0,\CA_0)$. Then the
index of $\hat\Gamma$ is smaller or equal to the index of
$\Gamma_0$:
\[i(\hat\Gamma)\leq i(\Gamma_0).\]
\end{prop}
\begin{proof} For any point $P$ in $\Omega$, by definition,
there exist at least two distinct partial isometries $a$ and $b$
in $\CA_0^{\pm 1}$ defined at $P$ and such that $P.a$ and $P.b$
also lie in $\Omega$. For any $n\in\N$, $\Omega$ is contained in
$\CF_n$: let $C_n$ be the connected component of $P$ in $\CF_n$.
There are at least two edges going out of the vertex $C_n$ of
$\Gamma_n$ corresponding to the restrictions of $a$ and $b$ to
$C_n$. This proves that the image of $\hat\Gamma$ by
$\hat\tau_n$ in $\Gamma_n$ does not contain vertices of valence
$0$ or $1$: $\hat\tau_n(\hat\Gamma)$ is a subgraph of the core
of $\Gamma_n$.

Let $\Theta_0$ be a finite set of vertices of $\hat\Gamma$. Let
$\Theta$ be a finite subgraph of $\hat\Gamma$ that contains
$\Theta_0$ and all edges incident to elements of $\Theta_0$. The
graph $\Theta$ exists because vertices of $\hat\Gamma$ have
finite valence bounded above by twice the cardinality of
$\CA_0$.

By Lemma~\ref{lem:inverselimit}, there exists $n\in\N$ such that
$\Theta$ is mapped injectively by $\hat\tau_n$ into $\Gamma_n$.
Arguing as in Lemma~\ref{lem:indexgraph} and using the
definition of the index of a graph, the following inequalities
hold:
\[
\sum_{x\in\Theta_0}
i_{\hat\Gamma}(x)=\sum_{x\in\Theta_0} i_\Theta(x)\leq
\sum_{x\in\Theta_0}i_{\core(\Gamma_n)}(\hat\tau_n(x))\leq
i(\core(\Gamma_n))=i(\Gamma_n).
\]

Using Proposition~\ref{prop:indexripsmachine} we get that for
any finite subset $\Theta_0$ of vertices of $\hat\Gamma$
\[
\sum_{x\in\Theta_0} i_{\hat\Gamma}(x)\leq i(\Gamma_0).
\]

Taking $\Theta_0$ arbitrarily large we finally get
\[
i(\hat\Gamma)\leq i(\Gamma_0). \qedhere
\]
\end{proof}

A connected component $K$ of $\Omega$ is called \textbf{regular}
if it corresponds to a vertex of valence $2$ of the limit graph
$\hat\Gamma$. A connected component $K$ of $\Omega$ is
\textbf{singular} if it corresponds to a vertex of valence at
least $3$ of $\hat\Gamma$.

\begin{cor}\label{cor:finitesingular} 
All connected components of $\hat\Gamma$ are lines except at most
$i(\Gamma_0)$.  Moreover, there are at most $i(\Gamma_0)$ singular
connected components of $\Omega$.\qed
\end{cor}

We are now ready to prove Theorem~\ref{thm:mainsi} in the
special case where the Rips Machine never halts and digs holes
everywhere.

\begin{thm}\label{thm:pseudo-levitt} 
Let $\CK=(\CF,\CA)$ be system of isometries and $\Gamma$ be its
associated graph.  Assume that the limit set $\Omega$ is totally
disconnected. Then the index of $\CK$ is bounded above by the index of
$\Gamma$.
\end{thm}
\begin{proof} As connected components of $\Omega$ are
reduced to single points, the graph $\hat\Gamma$ is the disjoint
union of all the cores of the Cayley graphs $\Cayl(P)$:
\[
\hat\Gamma=\biguplus_{[P]}\core(\Cayl(P)).
\]
Thus, the
index of $\hat\Gamma$ is equal to the index of $\CK$. The
Theorem now follows from Proposition~\ref{prop:hatgamma}.
\end{proof}

We turn back to the general case where the limit set has
non-trivial (i.e. not reduced to a single point) connected
components.

\begin{prop}\label{prop:finiteconnectedcomponents} The limit set
$\Omega$ of a system of isometries $\CK_0=(\CF_0,\CA_0)$ has
finitely many non-trivial connected components.
\end{prop}
\begin{proof} By Corollary~\ref{cor:finitesingular}, there are
finitely many singular connected components of $\Omega$.

Let $K$ be a regular connected component of $\Omega$. Then there
exists exactly two distinct partial isometries $a\neq b$ in
${\CA_0}^{\pm 1}$ with non trivial restrictions $a'=
{}_{K\rceil}a_{\lceil\Omega}$ and $b'=
{}_{K\rceil}b_{\lceil\Omega}$ to $K$ and $\Omega$. These are the
partial isometries which give rise to the two edges of
$\hat\Gamma$ going out of the vertex $K$. By definition of
$\Omega$, for each point $P$ of $K$ there exists at least two
partial isometries in $\CA_0^{\pm 1}$ which map $P$ inside
$\Omega$. Thus $P$ lies in both the domains of $a'$ and $b'$ and
thus $\dom (a')=\dom (b')=K$.

Now, if the range $K.a'$ of $a'$ is also a regular connected
component of $\Omega$ then $a'$ is an isometry between $K$ and
$K.a'$ and in particular they have the same diameter. From the
next Lemma we get that there can only be finitely many
non-trivial regular connected components.

This proves that $\Omega$ only has finitely many non-trivial
connected components.
\end{proof}

\begin{lem}\label{lem:finitediameter} 
Let $K$ be a compact $\R$-tree and $(K_i)_{i\in\N}$ be a collection of
disjoint subtrees of $K$. Then
\[
\lim_{i\to\infty}\mathrm{diam}(K_i)=0.
\]
\end{lem}
\begin{proof} By contradiction, assume that there
exists $\epsilon>0$ such that (extracting a subsequence)
$\forall i$, diam$(K_i)>\epsilon$. Let $P_i,Q_i$ be in $K_i$
with $d(P_i,Q_i)\geq \epsilon$. By compacity of $K$, up to
passing to a subsequence, the sequences $P_i$ and $Q_i$ converge
to $P$ and $Q$ respectively in $K$ and $d(P,Q)\geq\epsilon$. For
$i$ big enough $d(P_i,P)<\epsilon/2$ and $d(Q_i,Q)<\epsilon/2$,
from which we deduce that the segment $[P_i,Q_i]$ contains the
midpoint of $[P,Q]$, and therefore the segments $[P_i,Q_i]$ and
$[P_j,Q_j]$ are not disjoint for $i$ and $j$ big enough. A
contradiction.
\end{proof}

\subsection{Reduced systems of isometries}\label{sec:reduced}

In this section we introduce an extra hypothesis on the system
of isometries. Under this hypothesis using the Rips Machine is
considerably easier. For system of isometries on finite trees,
such an extra hypothesis was introduced by D.~Gaboriau in
\cite{gab-bouts} where it appears in Proposition~V.4 as
Property~(*).

The \textbf{set of extremal points} $\partial K$ of a compact
$\R$-tree $K$ is the set of points of $K$ that do not lie in the
interior of an interval contained in $K$. Equivalently a point
$P$ is extremal in $K$ if $K\smallsetminus\{P\}$ is connected.
The tree $K$ is the convex hull of its extremal points:
\[
K=\bigcup_{P,P'\in\partial K} [P;P'].
\]
We remark that
$\partial K$ may fail to be compact.

\begin{defn}\label{def:reduced} Let $\CK_0=(\CF_0,\CA_0)$ be a
system of isometries. Let $\CK_1=(\CF_1,\CA_1)$ be the output of
the Rips Machine. The system of isometries $\CK_0$ is
\textbf{reduced} if
\begin{enumerate}
\renewcommand{\theenumi}{(\roman{enumi})}
\renewcommand{\labelenumi}{\theenumi}
\item\label{defn:finitetrajectories} For any point $P$ in
$\CF_0$ the tree of trajectories $\Tt(P)$ is infinite.
\item\label{defn:convexlocal} For any partial isometry $a$ in
$A_0^{\pm 1}$ the set of extremal points $\partial\dom(a)$ of
the domain of $a$ is contained in $\CF_1$.
\end{enumerate}
\end{defn}

\begin{lem}\label{lem:reducedGamma} The graph $\Gamma_0$
associated to a reduced system of isometries does not have
vertices of negative index, that is to say of valence $0$ or
$1$.
\end{lem}
\begin{proof} Let $K$ be a connected component of
$\CF$ and let $P$ be a point in $K$. From
condition~\ref{defn:finitetrajectories}, the tree of
trajectories $\Tt(P)$ is infinite and thus there exists at least
one partial isometry $a\in\CA_0^{\pm 1}$ defined at $P$. Let now
$Q$ be an extremal point of $\dom(a)$. From
condition~\ref{defn:convexlocal}, $Q$ is in $\CF_1$ and thus
belongs to at least another domain $\dom(b)$ with
$b\in\CA_0^{\pm 1}$, $a\neq b$. Thus, the vertex $K$ of
$\Gamma_0$ has at least two outgoing edges $a$ and $b$.
\end{proof}

When dealing with the Rips Machine, reduced systems of isometries are
easier to handle. The elementary moves ``split a vertex'' and ``split
an edge'' are described in the proof of
Proposition~\ref{prop:indexripsmachine}.

\begin{prop}\label{prop:reducedRipsmachine} Let
$\CK_0=(\CF_0,\CA_0)$ be a reduced system of isometries. Let
$\CK_1=(\CF_1,\CA_1)$ be the output of the Rips Machine.

Then going from $\CK_0$ to $\CK_1$ only performs elementary
moves of type ``split an edge'', and the map
$\tau:\Gamma_1\to\Gamma_0$ is a homotopy equivalence. In
particular
\[
i(\Gamma_0)=i(\Gamma_1).
\]
\end{prop}
\begin{proof} In the proof of
Proposition~\ref{prop:indexripsmachine}, starting with a reduced
system of isometries, we first get that $E_0$ is empty. Indeed,
let $P$ be an extremal point of $\CF$ which is not in $\CF_1$.
From condition~\ref{defn:convexlocal}, $P$ does not belong to
any domain of a partial isometry in $\CA^{\pm 1}$ and thus its
tree of trajectories consists in a single vertex, which
contradicts condition~\ref{defn:finitetrajectories}.

Then we get that no ``split a vertex'' move can occur, because
the removed points in this move have a tree of trajectories
which consists in a single vertex.
\end{proof}

We can now prove that the reduced condition is inherited while
iterating the Rips Machine.

\begin{prop}\label{prop:reducedinherited} Let
$\CK_0=(\CF_0,\CA_0)$ be a reduced system of isometries. Let
$\CK_1=(\CF_1,\CA_1)$ be the output of the Rips Machine. Then
$\CK_1$ is reduced.
\end{prop}
\begin{proof} For a point $P$ in
$\CF_1$, the tree of trajectories of $P$ with respect to $\CK_1$
is obtained from the tree of trajectories with respect to
$\CK_0$ by pruning off the terminal vertices. Thus, if the
latter is infinite, so is the former.

By contradiction, let $\CK_2=(\CF_2,\CA_2)$ be the output of the
Rips Machine applied to $\CK_1$, let $a_1$ be a partial isometry
in $\CA_1$ and let $P$ be an extremal point in
$\partial\dom(a_1)$ which is not in $\CF_2$. Let $a_0$ be the
partial isometry in $\CA_0$ of which $a_1$ is a restriction.

As $P$ is in $\dom(a_1)$, and thus in $\CF_1$, there is at least
another partial isometry $b_0\in\CA_0$ such that $P$ is in the
domain $\dom(b_0)$. There exist two extremal points $Q,R$ in
$\partial\dom(b_0)$ such that $P$ lies in the segment $[Q;R]$.
By hypothesis, $Q,R$, as well as $Q.b_0$ and $R.b_0$ lie in
$\CF_1$.

As $P$ is an extremal point in $\partial\dom(a_1)$, up to
exchanging the two points $Q$ and $R$, we assume that $Q$ is not
in the domain of $a_1$ and that $[P;Q]\cap\dom(a_1)=\{P\}$. Let
$(Q_n)$ be a sequence of points in the open arc $]P;Q[$ that
converges to $P$. 

\begin{figure}[!h]
\centering
\setlength{\unitlength}{4144sp}%
\begingroup\makeatletter\ifx\SetFigFont\undefined%
\gdef\SetFigFont#1#2#3#4#5{%
  \reset@font\fontsize{#1}{#2pt}%
  \fontfamily{#3}\fontseries{#4}\fontshape{#5}%
  \selectfont}%
\fi\endgroup%
\begin{picture}(4094,2994)(429,-1873)
\thinlines
\put(901,839){\line( 0,-1){1170}}
\put(901,-331){\line( 1, 0){630}}
\put(1531,-331){\line( 0, 1){1170}}
\put(1801,-1681){\line( 0, 1){990}}
\put(1801,-691){\line( 1, 0){450}}
\put(2251,-691){\line( 0,-1){990}}
\put(3601,-1681){\line( 0, 1){990}}
\put(3601,-691){\line( 1, 0){450}}
\put(4051,-691){\line( 0,-1){990}}
\multiput(2971,-1141)(0.00000,-110.76923){7}{\line( 0,-1){ 55.385}}
\put(2971,-1861){\vector( 0,-1){0}}
\multiput(2251,389)(0.00000,114.54545){6}{\line( 0, 1){ 57.273}}
\put(2251,1019){\vector( 0, 1){0}}
\put(2071,-1231){\vector( 0,-1){450}}
\put(3871,-1141){\vector( 0,-1){540}}
\put(1261,119){\vector( 0, 1){630}}
\thicklines
\put(2881,-466){\line( 0,-1){ 90}}
\put(3331,-556){\line( 0, 1){ 90}}
\thinlines
\multiput(1711,-1861)(0.00000,109.56522){12}{\line( 0, 1){ 54.783}}
\multiput(1711,-601)(113.02326,0.00000){22}{\line( 1, 0){ 56.512}}
\multiput(4141,-601)(0.00000,-109.56522){12}{\line( 0,-1){ 54.783}}
\multiput(721,1109)(0.00000,-113.33333){14}{\line( 0,-1){ 56.667}}
\multiput(721,-421)(116.75676,0.00000){19}{\line( 1, 0){ 58.378}}
\multiput(2881,-421)(0.00000,113.33333){14}{\line( 0, 1){ 56.667}}
\thicklines
\put(451,-511){\line( 1, 0){4050}}
\put(1712,-467){\line( 0,-1){ 90}}
\put(4142,-472){\line( 0,-1){ 90}}
\put(2161,209){\makebox(0,0)[lb]{\smash{{\SetFigFont{12}{14.4}{\familydefault}{\mddefault}{\updefault}$a_1$}}}}
\put(1126,-61){\makebox(0,0)[lb]{\smash{{\SetFigFont{12}{14.4}{\familydefault}{\mddefault}{\updefault}$a_2$}}}}
\put(2881,-1006){\makebox(0,0)[lb]{\smash{{\SetFigFont{12}{14.4}{\familydefault}{\mddefault}{\updefault}$b_0$}}}}
\put(1981,-1096){\makebox(0,0)[lb]{\smash{{\SetFigFont{12}{14.4}{\familydefault}{\mddefault}{\updefault}$b_1$}}}}
\put(3691,-1006){\makebox(0,0)[lb]{\smash{{\SetFigFont{12}{14.4}{\familydefault}{\mddefault}{\updefault}$b'_1$}}}}
\put(3286,-421){\makebox(0,0)[lb]{\smash{{\SetFigFont{12}{14.4}{\familydefault}{\mddefault}{\updefault}$Q_n$}}}}
\put(2926,-421){\makebox(0,0)[lb]{\smash{{\SetFigFont{12}{14.4}{\familydefault}{\mddefault}{\updefault}$P$}}}}
\put(4096,-430){\makebox(0,0)[lb]{\smash{{\SetFigFont{12}{14.4}{\familydefault}{\mddefault}{\updefault}$Q$}}}}
\put(1669,-367){\makebox(0,0)[lb]{\smash{{\SetFigFont{12}{14.4}{\familydefault}{\mddefault}{\updefault}$R$}}}}
\end{picture}%
\end{figure}

In the original system of isometries $\CK_0$, $Q_n$ has an
infinite tree of trajectories, in particular there exists a
reduced admissible path $c_n.d_n$ of length $2$ in $\Gamma_0$
which is defined at $Q_n$. As $\Gamma_0$ is a finite graph, up
to passing to a subsequence we assume that this path is
constant: for any $n$ the reduced admissible path $c.d$ is
defined at $Q_n$. As domains of partial isometries are closed, the
point $P$ is in the domain of $c.d$. Since $c\inv\neq d$, by
definition of the Rips Machine, the point $Pc$ is in $\CF_1$ and
there is a restriction $c_1$ of $c$ in $A_1$ which is defined at
$P$. As $P$ is not in $\CF_2$ and as $a_1$ is defined at $P$,
the partial isometries $c_1$ and $a_1$ are equal. Thus for any
integer $n$ the partial isometry $a_0$ is defined at $Q_n$ and
$Q_na_0$ is in $\CF_1$. There are only finitely many
restrictions of $a_0$ to the connected components of $\CF_1$.
Thus $a_1$ is defined at some $Q_n$ and $P$ is not an extremal
point in $\partial\dom(a_1)$. A contradiction.
\end{proof}

We now state an equivalent characterization of reduced systems
of isometries.

\begin{prop}\label{prop:reducedlocal} A system of isometries
$\CK_0=(\CF_0,\CA_0)$ is reduced if and only if the two
following conditions hold:
\begin{enumerate}
\renewcommand{\theenumi}{(\roman{enumi})}
\renewcommand{\labelenumi}{\theenumi} 
\item For any point $P$ in
$\CF$ the tree of trajectories $\Tt(P)$ is infinite.
\renewcommand{\theenumi}{(\roman{enumi}')}
\item\label{defn:convex} For any partial isometry $a$ in $A^{\pm
1}$ the set of extremal points $\partial\dom(a)$ of the domain
of $a$ is contained in the limit set $\Omega_0$.
\end{enumerate}
\end{prop}
\begin{proof} As the limit set $\Omega_0$ of $\CK_0$
is contained in $\CF_1$ we get that condition~\ref{defn:convex}
implies condition~\ref{defn:convexlocal}.

Conversely, let $\CK_0=(\CF_0,\CA_0)$ be a reduced system of
isometries and let $\CK_n=(\CF_n,\CA_n)$ be the systems of
isometries obtained by iteratively applying the Rips machine.
From Proposition~\ref{prop:reducedinherited}, $\CK_n$ is
reduced. Let $P$ be an extremal point of $\dom(a_0)$ for a
partial isometry $a_0\in\CA_0^{\pm 1}$. By
condition~\ref{defn:convexlocal}, $P$ is in $\CF_1$ as well as
$P.a_0$, and thus there exists a restriction $a_1\in\CA_1^{\pm
1}$ of $a_0$ such that $P$ is an extremal point of $\dom(a_1)$.
By induction for each $n$ there exists a partial isometry $a_n$
in $\CA_n^{\pm 1}$ such that $P$ is an extremal point in
$\dom(a_n)$ and thus $P$ is in $\CF_n$. By
Lemma~\ref{lem:limitrips} we conclude that $P$ is in the limit
set $\Omega_0$.
\end{proof}

\section{Computing the index of a system of isometries}\label{sec:compindexsi}

\subsection{Systems of isometries on finite trees}

We translate and adapt in this section Proposition~6.1 of
\cite{glp-rips}.

A \textbf{finite tree} is an $\R$-tree which is the convex hull
of finitely many of its points. It has finitely many branch
points and is the metric realization of a simplicial finite
tree. A \textbf{finite forest} is a metric space with finitely
many connected components each of which is a finite tree.

We remark that removing the branch points, such a finite forest
gives a disjoint union of finitely many intervals. The integral
of a function defined on $\CF$ is the integral on this disjoint
union of intervals (with respect to Lebesgue measure).

Let $\CK=(\CF,\CA)$ be a system of isometries where $\CF$ is a
finite forest.

The \textbf{valence} of a point $P$ in $\CF$ is
\[
v_\CK(P)=\#\{a\in\CA^{\pm 1}\ |\ P\in\dom(a)\}.
\]
We emphasize
that in Section~\ref{sec:indexsi} we defined the index
$i_\CK(P)$ by restricting partial isometries to the limit set
$\Omega$, and thus for a system of isometries we may have
$i_\CK(P)\neq v_\CK(P)-2$.

The function $P\mapsto v_\CK(P)$ is a finite sum of
characteristic functions of finite subtrees. It is Lebesgue
integrable.

We translate Proposition~6.1 of \cite{glp-rips} to get:

\begin{prop}\label{prop:finitetree} Let $\CK=(\CF,\CA)$ be a
system of isometries with independent generators. Assume that
$\CF$ is a finite forest. Then $\int_{P\in \CF}(v_\CK(P)-2)\leq
0$
\end{prop}
\begin{proof} The sum $\int_{P\in\CF}v_\CK(P)$ is
twice the sum of the Lebesgue measure of the domains of the
partial isometries in $\CA$. It is denoted by $2\ell$ in
Section~6 of \cite{glp-rips} while $\int_{P\in\CF}2$ is twice
the Lebesgue measure of $\CF$ which is denoted there by $2m$.
\end{proof}

\subsection{Shortening systems of isometries}\label{sec:shortening}

To use Proposition~\ref{prop:finitetree} in the broader context
of compact forests, we need a standard procedure to shorten a
system of isometries to a finite one.

Let $\CF$ be a compact forest and let $\epsilon>0$. We define:
\[
\CF_{\chop\epsilon}=\{ P\in\CF\ |\ \exists P_{-\epsilon},
P_{+\epsilon}\in\CF, P\in [ P_{-\epsilon},P_{+\epsilon}],
d(P,P_{-\epsilon})=d(P,P_{+\epsilon})=\epsilon\}
\]
(i.e. $P$ is
the midpoint of the segment $[ P_{-\epsilon},P_{+\epsilon}]$
which has length $2\epsilon$).

The set of extremal points $\partial\CF$ of a compact forest
$\CF$ is the union of the sets of extremal points of the
connected components of $\CF$.

\begin{lem}\label{lem:chop}
\begin{enumerate} 
\item[(i)] For any
$\epsilon>0$, for any connected component $K$ of $\CF$,
$K\cap\CF_{\chop\epsilon}=K_{\chop\epsilon}$ is a finite tree
(possibly empty); 
\item[(ii)] For any $\epsilon>\epsilon'>0$,
$\CF_{\chop\epsilon}\subset \CF_{\chop\epsilon'}$; 
\item[(iii)]
$\cup_{\epsilon>0}\CF_{\chop\epsilon}=\CF\smallsetminus\partial\CF$.\qed
\end{enumerate}
\end{lem}

For any partial isometry $a$ of $\CF$, we denote by $a_\epsilon$
its restriction to $\CF_{\chop\epsilon}$. We denote by
$\CA_\epsilon$ and $\CK_\epsilon$ the corresponding finite set
of partial isometries and the restricted system of isometries.
We remark that if $\CK$ has independent generators then
$\CK_\epsilon$ also has.

\subsection{Pseudo-surface systems of isometries}\label{sec:pseudo-surface}

Let $\CK=(\CF,\CA)$ be a system of isometries, where $\CF$ is a
compact forest and such that each point of $\CF$ lies inside the
domain of at least two different partial isometries in $\CA^{\pm
1}$. In this case the limit set $\Omega$ is equal to $\CF$ and
the Rips Machine does not do anything to $\CK$.

If, in addition, the system of isometries $\CK$ has independent
generators, we say that it is \textbf{pseudo-surface}.

A key step in our proof of Theorem~\ref{thm:mainsi} is the
following Proposition which is proved using
Proposition~\ref{prop:finitetree} by shortening the system of
isometry as in Section~\ref{sec:shortening}. We note that the
following Proposition is obvious if $\CF$ is a finite tree or a
finite forest (cf. for example
\cite[Processus~II.3~5)]{gab-bouts}).

\begin{prop}\label{prop:thick} Let $\CK=(\CF,\CA)$ be a
pseudo-surface system of isometries. Then, for any choice of
three distinct partial isometries in $\CA^{\pm 1}$ the
intersection of the three domains contains at most one point of
$\CF$.
\end{prop}
\begin{proof} By definition of pseudo-surface
systems of isometries, for any $P$ in $\CF$ the valence $v_S(P)$
is greater or equal to $2$.

By contradiction we assume that there exist three distinct
partial isometries $a$, $b$ and $c$ in $\CA^{\pm 1}$ such that
the intersection of their domains is strictly bigger than a
point. As domains are compact subtrees, there exists a
non-trivial arc $I$ which is contained in the three domains:
\[
\forall P\in I, v_S(P)\geq 3.
\]

For any $\epsilon>0$ we consider the finite forest
$\CF_{\chop\epsilon}$ and the corresponding system of isometries
$S_\epsilon$.

Let $\ell$ be the length of $I$. For any $\epsilon<\ell/3$, the
sub-arc $J$ of $I$ which is contained in the domains of the
partial isometries $a_\epsilon$, $b_\epsilon$ and $c_\epsilon$
of $\CA_\epsilon$ contains the middle third of $I$ thus the
length of $J$ is bigger than $\ell/3$.

For such an $\epsilon>0$, $\CF_{\chop\epsilon}$ is a finite
forest by Lemma~\ref{lem:chop} and
Proposition~\ref{prop:finitetree} holds:
\[
I_\epsilon=\int_{P\in\CF_{\chop\epsilon}}(v_{\CK_\epsilon}(P)-2)
\leq 0
\]
Let $\CP_\epsilon$, resp. $\CN_\epsilon$, be the set
of points of $\CF_{\chop\epsilon}$ which contributes positively,
resp. negatively, to $I_\epsilon$:
\[
\CP_\epsilon=\{
P\in\CF_{\chop\epsilon}\ |\ v_{\CK_\epsilon}(P)\geq 2\}\mbox{
and }\CN_\epsilon=\{ P\in\CF_{\chop\epsilon}\ |\
v_{\CK_\epsilon}(P)\leq 1\}.
\]
As the points in $J$ have
valence at least $3$, we get
\[
\int_{P\in\CP_\epsilon}(v_{\CK_\epsilon}(P)-2)\geq \ell/3
\]
and
thus, we have:
\[
0\geq I_\epsilon\geq
\ell/3+\int_{P\in\CN_\epsilon}(v_{\CK_\epsilon}(P)-2).
\]
Our
goal is to prove that the negative part goes to zero to get a
contradiction. We only need to prove that the Lebesgue measure
of $\CN_\epsilon$ goes to zero.

We claim that $\CN_\epsilon$ has Lebesgue measure bounded above
by $8N^2\,\epsilon$ where $N=\#\CA$ is the cardinality of $\CA$.

Let $P$ be in $\CN_\epsilon$, then $P$ is in
$\CF_{\chop\epsilon}$ which means that $P$ is the midpoint of a
segment $[P_{-\epsilon}, P_{+\epsilon}]$ of length $2\epsilon$
in $\CF$. As $S$ is pseudo-surface, there are at least two
elements $a_1,a_2\in \CA^{\pm 1}$ which are defined at
$P_{-\epsilon}$ and at least two elements $b_1,b_2\in \CA^{\pm
1}$ which are defined at $P_{+\epsilon}$. As $P$ is in
$\CN_\epsilon$ at most one of the four partial isometries
$a_{1\epsilon}$,$a_{2\epsilon}$,$b_{1\epsilon}$,$b_{2\epsilon}$
is defined at $P$. By switching the indices we can assume that
$a_{1\epsilon}$ and $b_{1\epsilon}$ are not defined at $P$.
Taking $\epsilon$ sufficiently small ensures that the partial
isometries $a_{1\epsilon}$ and $b_{1\epsilon}$ are not empty.

The domain of $a_{1\epsilon}$ and the point $P_{-\epsilon}$ lie
in the same connected component of $F\smallsetminus\{P\}$: else
$P$ would be located in a segment $[P_{-\epsilon},P']$ with $P'$
in the domain of $a_{1\epsilon}$ and $P_{-\epsilon}$ in the
domain of $a_1$, thus $P$ would be in the domain of
$a_{1\epsilon}$. We argue similarly for the domain of
$b_{1\epsilon}$. 

\begin{figure}[h]
\centering
\setlength{\unitlength}{4144sp}%
\begingroup\makeatletter\ifx\SetFigFont\undefined%
\gdef\SetFigFont#1#2#3#4#5{%
  \reset@font\fontsize{#1}{#2pt}%
  \fontfamily{#3}\fontseries{#4}\fontshape{#5}%
  \selectfont}%
\fi\endgroup%
\begin{picture}(5466,1800)(643,-2299)
\thinlines
\put(2026,-1231){\vector( 0, 1){450}}
\put(5401,-1231){\vector( 0, 1){450}}
\put(1351,-1771){\line( 0,-1){180}}
\put(2926,-1771){\line( 0,-1){180}}
\put(4501,-1771){\line( 0,-1){180}}
\thicklines
\put(676,-1861){\line( 1, 0){5400}}
\thinlines
\multiput(901,-511)(0.00000,-109.56522){12}{\line( 0,-1){ 54.783}}
\multiput(901,-1771)(116.12903,0.00000){16}{\line( 1, 0){ 58.065}}
\multiput(2701,-1771)(0.00000,109.56522){12}{\line( 0, 1){ 54.783}}
\multiput(4321,-511)(0.00000,-109.56522){12}{\line( 0,-1){ 54.783}}
\multiput(4321,-1771)(111.72414,0.00000){15}{\line( 1, 0){ 55.862}}
\multiput(5941,-1771)(0.00000,109.56522){12}{\line( 0, 1){ 54.783}}
\multiput(1351,-961)(0.00000,128.57143){4}{\line( 0, 1){ 64.286}}
\put(1351,-511){\vector( 0, 1){0}}
\multiput(4771,-961)(0.00000,128.57143){4}{\line( 0, 1){ 64.286}}
\put(4771,-511){\vector( 0, 1){0}}
\put(5079,-736){\line( 0,-1){990}}
\put(5079,-1726){\line( 1, 0){630}}
\put(5709,-1726){\line( 0, 1){990}}
\put(1707,-736){\line( 0,-1){990}}
\put(1707,-1726){\line( 1, 0){765}}
\put(2472,-1726){\line( 0, 1){990}}
\put(1261,-2221){\makebox(0,0)[lb]{\smash{{\SetFigFont{12}{14.4}{\familydefault}{\mddefault}{\updefault}$P_{-\epsilon}$}}}}
\put(4366,-2221){\makebox(0,0)[lb]{\smash{{\SetFigFont{12}{14.4}{\familydefault}{\mddefault}{\updefault}$P_{+\epsilon}$}}}}
\put(2791,-2221){\makebox(0,0)[lb]{\smash{{\SetFigFont{12}{14.4}{\familydefault}{\mddefault}{\updefault}$P$}}}}
\put(1891,-1501){\makebox(0,0)[lb]{\smash{{\SetFigFont{12}{14.4}{\familydefault}{\mddefault}{\updefault}$a_{1\epsilon}$}}}}
\put(5221,-1501){\makebox(0,0)[lb]{\smash{{\SetFigFont{12}{14.4}{\familydefault}{\mddefault}{\updefault}$b_{1\epsilon}$}}}}
\put(1171,-1141){\makebox(0,0)[lb]{\smash{{\SetFigFont{12}{14.4}{\familydefault}{\mddefault}{\updefault}$a_1$}}}}
\put(4591,-1141){\makebox(0,0)[lb]{\smash{{\SetFigFont{12}{14.4}{\familydefault}{\mddefault}{\updefault}$b_1$}}}}
\end{picture}%
\end{figure}

We have thus proved that $P$ is in the non-trivial arc joining
the disjoint domains of $a_{1\epsilon}$ and $b_{1\epsilon}$. The
point $P_{-\epsilon}$ is in the domain of $a_1$ and thus at
distance less than $\epsilon$ of the domain of $a_{1\epsilon}$.
Thus $P$ is at distance less than $2\epsilon$ from this domain.
Hence, the length of the arc joining the disjoint domains of
$a_{1\epsilon}$ and $b_{1\epsilon}$ is at most $4\epsilon$.

If, over all the possible pairs of partial isometries, we sum
the lengths of the arcs, we get that the volume of
$\CN_\epsilon$ is bounded above by $(2N(2N-1)/2)\times
4\epsilon$. Which proves the claim and concludes the proof.
\end{proof}

If the compact forest $\CF=I$ is an interval and if the system
of isometries $\CK=(\CF,\CA)$ is pseudo-surface,
Proposition~\ref{prop:thick} states that this is the classical
case of an interval exchange transformation and $\CK$ is usually
called surface. This justifies the terminology of pseudo-surface
system of isometries.

From Proposition~\ref{prop:thick}, it is easy to deduce a rough
bound of the index of a pseudo-surface system of isometries.

\begin{cor}\label{cor:pseudo-surfacerough} Let $\CK=(\CF,\CA)$
be a pseudo-surface system of isometries. Then $i(\CK)$ is
finite and bounded above by a constant depending only on the
cardinality of $\CA$.
\end{cor}
\begin{proof} We denote, as
before, by $N=\#\CA$ the cardinal of $\CA$. From the previous
Proposition, there are at most $\binom{2N}{3}$ points in $\CF$
which belongs to the domains of at least three different partial
isometries in $\CA^{\pm 1}$. Each of these points has valence at
most $2N$. Adding up we get that
\[
i(\CK)=\sum_{P\in\CF}i_\CK(P)=\sum_{P\in\CF}(v_\CK(P)-2)\leq
N(2N-1)(2N-2)^2/3.\qedhere
\]
\end{proof}

We now state a combinatorial Lemma. 

Let $K$ be a compact $\R$-tree and let ${\mathcal
K}=(K_a)_{a\in\CA}$ be a finite collection of compact subtrees
of $K$. For such a collection, and for any point $P\in K$ we
denote by $v_{\mathcal K}(P)$ the number of elements $a$ of
$\CA$ such that $P$ is in $K_a$.

\begin{lem}\label{lem:indexsurface} Let ${\mathcal
K}=(K_a)_{a\in\CA}$ be a finite collection of compact subtrees
of a compact $\R$-tree $K$. Assume that
\begin{enumerate} \item
for any choice of three distinct elements of $\CA$ the
intersection of the corresponding subtrees is at most one point,
\item any element $P$ of $K$ is in at least two compact subtrees
$K_a$ and $K_b$ ($a\neq b\in\CA$).
\end{enumerate} Then
\[
\sum_{P\in K}(v_{\mathcal K}(P)-2)=\#\CA-2.
\]
\end{lem}
\begin{proof} Let $T$ be the convex hull in $K$ of all elements
$P\in K$ such that $v_{\mathcal K}(P)\geq 3$. From our first
hypothesis $T$ is a finite tree. For each $a\in\CA$ the
intersection $T_a=K_a\cap T$ is a finite tree. Let ${\mathcal
T}=(T_a)_{a\in\CA}$ be the corresponding collection of finite
subtrees of $T$. We have the equality
\[
\sum_{P\in
K}(v_{\mathcal K}(P)-2)=\sum_{P\in T}(v_{\mathcal T}(P)-2).
\]
Moreover $\mathcal T$ satisfies the same hypothesis as $\mathcal
K$: the intersection of three of its elements is at most a point
and any element of $T$ is in at least two subtrees $T_a$ and
$T_b$, $a\neq b\in\CA$. We regard $T$ and each $T_a$ as a
simplicial tree by considering all the branch points and
extremal points as vertices. Each edge of the simplicial tree
$T$ belongs to exactly two trees $T_a$ and $T_b$. Combinatorial
computation gives
\begin{eqnarray*} \sum_{P\in T}(v_{\mathcal
T}(P)-2)&=&\sum_{P\in V(T)}(v_{\mathcal T}(P)-2)\\ &=&\sum_{P\in
V(T)}\sum_{a\in A} \carac_{P\in V(T_a)}-2\#V(T)\\ &=&\sum_{a\in
A}\sum_{P\in V(T)} \carac_{P\in V(T_a)} - 2\#E(T)-2\\
&=&\sum_{a\in A}\#V(T_a)-\sum_{a\in A} \#E(T_a)-2\\ &=&\#\CA-2.
\end{eqnarray*}
\qedhere
\end{proof}

We can now get the correct bound for the index of a
pseudo-surface system of isometries

\begin{thm}\label{thm:thick} Let $\CK=(\CF,\CA)$ be a
pseudo-surface system of isometries and, let $\Gamma$ be the
associated graph. Then $i(\CK)=i(\Gamma)$.
\end{thm}
\begin{proof} As $\CK$ is pseudo-surface, at least two distinct
partial isometries are defined at each point of $\CF$. Thus, the
graph $\Gamma$ does not have vertices of valence $0$ or $1$ and
its index is given by
\[
i(\Gamma)=\sum_{K\in
V(\Gamma)}i_\Gamma(K).
\]

Let $K$ be a connected component of $\CF$ and let $B$ be the
subset of $A^{\pm 1}$ which consists of partial isometries with
domains inside $K$. The set $B$ is also the set of edges going
out of the vertex $K$ of the graph $\Gamma$ and thus
\[
i_\Gamma(K)=\#B-2.
\]
Let $\mathcal K$ be the collection of
domains of elements of $B$. Thus for each point $P$ in $K$
\[
i_S(P)=v_{\mathcal K}(P)-2.
\]
By Proposition~\ref{prop:thick},
the collection $\mathcal K$ satisfies the hypothesis of
Lemma~\ref{lem:indexsurface} and we get
\[
\sum_{P\in
K}(v_{\mathcal K}(P)-2)=\#B-2.
\]

Thus the contribution of the points of $K$ to the index of $\CK$
is equal to the contribution of the corresponding vertex of
$\Gamma$:
\[
\sum_{P\in K}i_S(P)=i_\Gamma(K).
\]
Adding up for
all connected components $K$ of $\CF$, proves the Theorem.
\end{proof}

\subsection{Proof of Theorem~\ref{thm:mainsi}}

Using Proposition~\ref{prop:hatgamma}, Theorem~\ref{thm:mainsi}
is a consequence of

\begin{prop}\label{prop:indexhatgamma} Let $\CK=(\CF,\CA)$ be a
system of isometries with independent generators. Let
$\hat\Gamma$ be its limit graph. Then the index of $\CK$ is
equal to the index of $\hat\Gamma$:
\[
i(\CK)=i(\hat\Gamma).
\]
\end{prop}
\begin{proof} Let $\Omega_0$ be the union of all
non-trivial connected components of the limit set $\Omega$. By
Proposition~\ref{prop:finiteconnectedcomponents}, $\Omega_0$ has
finitely many connected components, that is to say, $\Omega_0$
is a compact forest. Let $\CK_0=(\Omega_0,\CA_0)$ be the system
of isometries which consists of the restrictions of $\CK$ to
$\Omega_0$. By definition of the limit set, $\CK_0$ is a
pseudo-surface system of isometries. By
Proposition~\ref{prop:thick}, the intersection the domains of
three distinct partial isometries of ${\CA_0}^{\pm 1}$ contains
at most one point.

For a vertex of $\hat\Gamma$ corresponding to a connected
component $K$ (possibly a single point) of $\Omega$ we can apply
Lemma~\ref{lem:indexsurface} to the collection given by the
domains of the edges going out of $K$ to get:
\[
i_{\hat\Gamma}(K)=\sum_{P\in K}(v_\CK(P)-2).
\]
By
Corollary~\ref{cor:finitesingular}, $\hat\Gamma$ has finitely
many vertices with index strictly positive (and these indices
are finite). Adding up for all these singular vertices of
$\hat\Gamma$, we get
\[
i(\CK)=i(\hat\Gamma).\qedhere
\]
\end{proof}

\subsection{Systems of isometries with maximal
index}\label{sec:simaxindex}

From Theorem~\ref{thm:mainsi}, we say that a system of
isometries $\CK$ has \textbf{maximal index} if its index is
equal to the index of its associated graph $\Gamma$:
$i(\CK)=i(\Gamma)$.

The following Proposition characterizes reduced systems of
isometries with maximal index.

\begin{thm}\label{thm:maxindex} Let $\CK=(\CF,\CA)$ be a reduced system
  of isometries with independent generators, let $\Gamma$ be its
  associated graph, and $\hat\Gamma$ be its limit graph. The following
  are equivalent
\begin{enumerate} \item\label{item:indexmax} The system of
isometries $\CK$ has maximal index,
\item\label{item:hatgammafinite} The graph $\hat\Gamma$ is
finite, \item\label{item:ripsfinite} The Rips Machine, starting
from $\CK$, halts after finitely many steps.
\end{enumerate}
\end{thm}
\begin{proof} As before we denote by
$\CK_n=(\CF_n,\CA_n)$ the system of isometries obtained after
$n$ steps of the Rips Machine. By
Proposition~\ref{prop:reducedRipsmachine} and
Proposition~\ref{prop:reducedinherited} the Rips Machine only
performs ``split an edge'' moves and induces a homotopy
equivalence $\tau_n:\Gamma_{n+1}\to\Gamma_n$ at each step. And
thus the index is constant:
\[
i(\Gamma_n)=i(\Gamma_0)=i(\Gamma).
\]
Moreover at each step the
graph $\Gamma_n$ does not have vertices of valence $0$ or $1$.

\underline{\ref{item:ripsfinite}$\Rightarrow$\ref{item:hatgammafinite}:}
If the Rips Machine halts after step $n$, then
$\Gamma_{n+1}=\Gamma_n=\hat\Gamma$ is a finite graph.

\underline{\ref{item:hatgammafinite}$\Rightarrow$\ref{item:ripsfinite}:}
Conversely, at each step $n$ the Rips Machine only performs
``split an edge'' moves. This move adds one edge to $\Gamma_n$
to get $\Gamma_{n+1}$. If the Rips Machine never halts, then the
number of edges of $\Gamma_n$ goes to infinity. As each of the
$\tau_n$ is onto, we get that $\hat\Gamma$ is infinite.

\underline{\ref{item:ripsfinite}$\Rightarrow$\ref{item:indexmax}:}
If the Rips Machine halts after finitely many
steps: for $n$ big enough
\[\hat\Gamma=\Gamma_{n+1}=\Gamma_n,\]
by Proposition~\ref{prop:indexhatgamma}, $i(\hat\Gamma)=i(\CK)$
and we get that $\CK$ has maximal index.

\underline{\ref{item:indexmax}$\Rightarrow$\ref{item:ripsfinite}:}
Assume that $\CK$ has maximal index:
$i(S)=i(\hat\Gamma)=i(\Gamma)$.

We proceed as in the proof of Proposition~\ref{prop:hatgamma}.
Let $\Theta_0$ be the finite subset of vertices of $\hat\Gamma$
with valence strictly bigger than $2$. Let $\Theta$ be the
finite subgraph of $\hat\Gamma$ which contains all edges
incident to $\Theta_0$. The graph $\Theta$ contains all the
index of $\hat\Gamma$:
\[
i(\hat\Gamma)=\sum_{x\in
\Theta_0}i_\Theta(x).
\]
For $n$ big enough, $\hat\tau_n$ is
injective on $\Theta$ and thus for each vertex $x$ in $\Theta_0$
\[
i_\Theta(x)\leq i_{\Gamma_n}(\hat\tau_n(x)).
\]
We assumed
that $\CK$ is reduced and thus that $\Gamma_n$ does not have
vertices of strictly negative index.

By maximality of the index, $i(\hat\Gamma)=i(\Gamma_n)$ and thus
we can compute
\[
i(\hat\Gamma)=\sum_{x\in\Theta_0}i_{\Theta}(x)=\sum_{x\in\Theta_0}i_{\Gamma_n}(\hat\tau_n(x))+\sum_{y\in
V(\Gamma_n)\smallsetminus
\hat\tau_n(\Theta_0)}i_{\Gamma_n}(y)=i(\Gamma_n).
\]
We deduce
that for each $x$ in $\Theta_0$ and for all $y\in
V(\Gamma_n)\smallsetminus \hat\tau_n(\Theta_0)$,
\[
i_\Theta(x)=i_{\Gamma_n}(\hat\tau_n(x))\mbox{ and
}i_{\Gamma_n}(y)=0.
\]
The commutative diagram
\[
\xymatrix{
&\hat\Gamma\ar[ld]_{\hat\tau_{n+1}} \ar[rd]^{\hat\tau_n}\\
\Gamma_{n+1}\ar[rr]^{\tau_n}&&\Gamma_n\\ }
\]
restricts to graph
isomorphisms between $\Theta$ and its images. Moreover,
$\hat\tau_{n+1}(\Theta_0)$ and $\hat\tau_n(\Theta_0)$ contain
all the vertices of strictly positive index of $\Gamma_{n+1}$
and $\Gamma_n$ respectively. Thus no ``split an edge'' move can
occur when passing from $\Gamma_n$ to $\Gamma_{n+1}$ and thus
the Rips Machine does not do anything to the system of
isometries $\CK_n$.
\end{proof}

\section{Trees}

Throughout this Section, $T$ is an $\R$-tree with a very small,
minimal action of the free group $\FN$ of rank $N$ by isometries
with dense orbits.

\subsection{The map $\CQ$}\label{sec:Q}

Let $P$ be a point in $T$, we consider the equivariant map
$\CQ_P:\FN\to T$, $u\mapsto u\cdot P$. This maps does not
extends continuously to the boundary $\partial\FN$ of $\FN$. To
overcome this difficulty we weaken the topology on $T$ by
considering the observers' topology.

Let $\hat T=\bar T\cup\partial T$ be the union of the metric
completion of $T$ and its (Gromov) boundary. $\hat T$ inherits
from the metric on $T$ a well defined topology. However, $\hat
T$ is not compact in general.

We consider on $\hat T$ the weaker \textbf{observers' topology}
and we denote by $\Tobs$ this topological space. A basis of open
sets in $\Tobs$ is given by the set of connected components of
$\hat T\smallsetminus\{P\}$ for all points $P$. This topology is
Hausdorff and $\Tobs$ is a compact space with the same connected
subspaces than $\hat T$, see \cite{chl2}.

\begin{prop}[\cite{chl2}]\label{prop:Qexists} Let $T$ be an
$\R$-tree with a very small, minimal action of $\FN$ by
isometries with dense orbits. There exists a unique map $\CQ$
that is the continuous extension from $\partial\FN$ to $\Tobs$
of the map $\CQ_P:u\mapsto u\cdot P$. The map $\CQ$ does not
depend on the choice of a point $P$. \qed
\end{prop}

This map $\CQ$ was first introduced by Levitt and Lustig in
\cite{ll-north-south, ll-periodic} with a slightly different
approach. In particular they proved

\begin{prop}\label{prop:Qonto} Let $T$ be an $\R$-tree with a
very small, minimal action of $\FN$ by isometries with dense
orbits. The map $\CQ$ is onto $\hat T$. The points $P$ in $\hat
T$ with strictly more than one pre-image by $\CQ$ are in the
metric completion $\bar T$ of $T$ (and not in the (Gromov)
boundary $\partial T$).
\end{prop}

It has been asked by Levitt and Lustig
\cite[Remark~3.6]{ll-north-south} whether the map
$\CQ:\partial\FN\to\hat T$ has finite fibers (in the case where
the action is free). We are going to answer this question and to
give a precise bound for the cardinal of the fibers. In this
purpose we need to make this question precise by the following
definition of the $\CQ$-index.

\subsection{The $\CQ$-index}\label{subsec:Q-index}

We denote by $\Stab(P)$ the stabilizer in $\FN$ of a point $P$
of $\hat T$.

It is proved in \cite{gl-rank} that $\Stab(P)$ is a finitely
generated subgroup of $\FN$. The subgroup $\Stab(P)$ is a free
group and its boundary $\partial\Stab(P)$ embeds in the boundary
of $\FN$. For any element
$X\in\partial\Stab(P)\subseteq\partial\FN$,
Proposition~\ref{prop:Qexists} proves that $\CQ(X)=P$. Elements
of $\partial\Stab(P)$ are called \textbf{singular}, and the
other elements of the fiber $\CQ\inv(P)$ are \textbf{regular}.
We denote by $\CQ\inv_r(P)$ the set of regular points. As $\CQ$
is equivariant, $\Stab(P)$ acts on $\CQ\inv(P)$ and on
$\CQ\inv_r(P)$.

The \textbf{$\CQ$-index} $\iq(P)$ of a point $P$ in $\hat T$
is defined by
\[
\iq(P)=\#(\CQ\inv_r
(P)/\Stab(P))+2\,\rank(\Stab(P))-2.
\]
When $\Stab(P)$ is
trivial this definition becomes
\[
\iq(P)=\#\CQ\inv (P)-2.
\]

The $\CQ$-index only depends on the orbit $[P]$ of $P$ under the
action of $\FN$ and we can define the \textbf{$\CQ$-index} of
the tree $T$ by
\[
\iq(T)=\sum_{[P]\in \hat T/\FN}\max
(0;i_{\CQ}([P])).
\]
From Proposition~\ref{prop:Qonto}, points in
$\partial T$ have exactly one pre-image by $\CQ$. Thus, only
points in $\bar T$ contribute to the $\CQ$-index of $T$.

The main goal of this section is to prove the following Theorem:
\begin{thm}\label{thm:iq} Let $T$ be an $\R$-tree with a very
small, minimal action of $\FN$ by isometries with dense orbits.
Then
\[
\iq(T)\leq 2N-2.
\]
\end{thm}
 
In the case of a free action of the free group $\FN$ on $T$ this
gives the answer to Levitt and Lustig's question:

\begin{cor}\label{cor:Qindexfinite} Let $T$ be an $\R$-tree with
a free, minimal action of $\FN$ by isometries with dense orbits.
Then, there are finitely many orbits of points $P$ in $\hat T$
with strictly more than $2$ elements in their $\CQ$-fiber
$\CQ\inv(P)$ and these fibers are finite.\qed
\end{cor}

\subsection{Dual lamination and compact heart}\label{subsec:Q2}

The \textbf{double boundary} of $\FN$ is
\[
\partial^2\FN=(\partial\FN\times\partial\FN)\smallsetminus\Delta
\]
where $\Delta$ is the diagonal. An element of $\partial^2\FN$
is a \textbf{line}.

Using the map $\CQ$, in \cite{chl1-II}, the \textbf{dual
lamination} $L(T)$ to the tree $T$ is defined.
\[
L(T)=\{(X,Y)\in\partial^2\FN\ |\ \CQ(X)=\CQ(Y)\}.
\]
From this definition, the map $\CQ$ naturally induces an equivariant
map $\CQ^2: L(T)\to \hat T$. It is proved in \cite{chl1-II} that the
map $\CQ^2$ is continuous and its image is a subset $\Omega$ of $\bar
T$ which we call the \textbf{limit set}.

We fix a basis $\CA$ of $\FN$. Elements of $\FN$ are reduced
finite words in the alphabet $\CA^{\pm 1}$. An element $X$ of
$\partial\FN$ is an infinite reduced words in $\CA^{\pm 1}$, we
denote by $X_1$ its first letter.

The \textbf{unit cylinder} $C_\CA(1)$ of $\partial^2\FN$ is
\[
C_\CA(1)=\{(X,Y)\in\partial^2\FN\ |\ X_1\neq Y_1\}.
\]
Although
$\partial^2\FN$ is not compact, the unit cylinder is compact and
indeed a Cantor set.

In \cite{chl4} the \textbf{relative limit set} of $T$ with
respect to $\CA$ is defined:
\[
\Omega_\CA=\CQ^2(L(T)\cap
C_\CA(1)).
\]
From the continuity of the map $\CQ^2$, the
relative limit set $\Omega_\CA$ is a compact subset of $\bar T$.

The \textbf{compact heart} $K_\CA$ of $T$ is the convex hull of
$\Omega_\CA$.

For any element $a$ of the basis $\CA$ we consider the partial
isometry (which we also denote by $a$, but which we write on the
right) which is the restriction of the action of $a\inv$:
\[
\begin{array}{rcl} K_\CA\cap aK_\CA&\to&K_\CA\cap a\inv K_\CA\\
x&\mapsto&x.a=a\inv x \end{array}
\]
We get a system of
isometries $\CK_\CA=(K_\CA,\CA)$ as defined in
Section~\ref{sec:si}.

In \cite{chl4} it is proved that $\CK$ encodes all the
informations given by $T$ and the action of $\FN$. To be more
specific, we summarize results of \cite{chl4} as follows:

\begin{prop}[\cite{chl4}]\label{prop:chl4} Let $T$ be an
$\R$-tree with a very small, minimal action of $\FN$ by
isometries with dense orbits. Let $\CA$ be a basis of $\FN$, let
$K_\CA$ be its compact heart and $\CK_\CA=(K_\CA,\CA)$ be the
associated system of isometries. Then
\begin{enumerate}
\item\label{item:chl4independent} $\CK_\CA$ has independent
generators. \item\label{item:chl4admissible} For any word
$u\in\FN$, and for any point $P\in K_\CA$, $u\inv.P\in K_\CA$ if
and only if $u$ is admissible for $\CK_\CA$ and $P\in\dom(u)$.
In this case $P.u=u\inv.P$. \item\label{item:chl4Q} For any
element $X\in\partial\FN$, $\CQ(X)=P\in K_\CA$ if and only if
$X$ is admissible and $\{P\}=\dom(X)$.
\end{enumerate}
\end{prop}
\begin{proof} Assertion~\ref{item:chl4independent} is
Lemma~5.1 of \cite{chl4}. Assertion~\ref{item:chl4admissible} is
proved in Lemma~3.5 (1) and Corollary~5.5 and
Assertion~\ref{item:chl4Q} is a consequence of Proposition~4.3,
Lemma~4.7 and Corollary~5.5.
\end{proof}

We deduce that for an infinite reduced admissible word $X$ the
definition of $\CQ(X)$ of Section~\ref{subsec:sidef} agrees with
the definition given by Proposition~\ref{prop:Qexists}. Moreover
the relative limit set $\Omega_\CA$ of the $\R$-tree $T$ is
equal to the limit set of the system of isometries $\CK_\CA$.

\subsection{The compact heart of a tree is
reduced}\label{sec:heartreduced}

As explained in Section~\ref{sec:reduced}, reduced systems of
isometries (see Definition~\ref{def:reduced}) considerably
simplifies the use of the Rips Machine. Fortunately, in the
context of $\R$-trees, which we are studying, and thanks to
\cite{chl4}, we can work with reduced systems of isometries.

\begin{prop}\label{prop:reducedheart} Let $T$ be an $\R$-tree
with a minimal very small action of $\FN$ by isometries with
dense orbits. Let $\CA$ be a basis of $\FN$, let $\Omega_\CA$ be
the relative limit set and $K_\CA$ be the compact heart of $T$.
Let $\CK_\CA=(K_\CA,\CA)$ the induced system of isometries.

Then, the system of isometries $\CK_\CA$ is reduced.
\end{prop}
\begin{proof} By Proposition~\ref{prop:Qonto}, the map $\CQ$ is
onto $\hat T$: for any point $P$ in $K_\CA$ there exists
$X\in\partial\FN$ such that $\CQ(X)=P$. By
Proposition~\ref{prop:chl4}, $X$ is admissible and $P$ is in the
domain of $X$. Thus, any point $P$ of the compact heart $K_\CA$
has an infinite tree of trajectories.
  
Let $\CK_1=(\CF_1,\CA_1)$ be the output of the Rips Machine. Let
$a$ be a partial isometry in $\CA$ and let $P$ be an extremal
point of the domain of $a$. If $P$ is an extremal point of
$K_\CA$, then, as $K_\CA$ is the convex hull of the relative
limit set $\Omega_\CA$, we get that $P$ is in $\Omega_\CA$ and
thus in $\CF_1$. If $P$ is not an extremal point of $K_\CA$
there exists a sequence $Q_n$ of points in $\CF$ which converges
to $P$ and which are not in the domain of $a$. The points in the
compact heart $K_\CA$ have infinite tree of trajectories, thus
for each of these points $Q_n$ there exists a partial isometry
$b_n$ in $\CA$ such that $Q_n$ is in the domain of $b_n$. As
$\CA$ is finite, up to passing to a subsequence we assume that
all the $Q_n$ are in the domain of a partial isometry $b$ in
$\CA$. The domain of $b$ is close and thus $P$ is in the domain
of $b$. By definition of the Rips Machine we get that $P$ is in
$\CF_1$.
\end{proof}

\subsection{Proof of Theorem~\ref{thm:iq}}

Theorem~\ref{thm:iq} is a consequence of
Theorem~\ref{thm:mainsi} and of the following Theorem which
relates the $\CQ$-index of $T$ and the index of the system of
isometries $\CK_\CA$ defined on its compact heart.

\begin{thm}\label{thm:siiq} Let $T$ be an $\R$-tree with a very
small, minimal action of $\FN$ by isometries with dense orbits.
Let $\CA$ be a basis of $\FN$. The $\CQ$-index of $T$ and the
index of the induced system of isometries $\CK_\CA=(K_\CA,\CA)$
on the heart $K_\CA$ of $T$ and $\CA$ are equal:
\[
\iq(T)=i(\CK_\CA).
\]
\end{thm}
\begin{proof} Let $P$ be a
point in $\hat T$ and $[P]$ be its orbit under the action of
$\FN$. By Proposition~\ref{prop:chl4}
assertion~\ref{item:chl4admissible}, the intersection of the
orbit $[P]$ and of the compact heart $K_\CA$ is a pseudo-orbit
(possibly empty) of the system of isometries $\CK_\CA$.

Assume that $\iq([P])\geq 0$. There are at least two distinct
elements $X,Y$ in the fiber $\CQ\inv(P)$. Let $u$ be the common
prefix of $X$ and $Y$, then $X'=u\inv X$ and $Y'=u\inv Y$ have
different first letter and are in the pre-image by $\CQ$ of
$P'=u\inv P$. By definition $P'$ is in the relative limit set
$\Omega_\CA$ of the tree $T$.

This proves that the $\CQ$-index of $T$ can be computed by
considering only pseudo-orbits in $K_\CA$:
\[
\iq(T)=\sum_{[P]\in K_\CA/\FN}\max(0;\iq([P])).
\]

Let $P$ be a point in $K_\CA$. By Proposition~\ref{prop:chl4},
the boundary at infinity of the tree of trajectories $\Tt(P)$ is
exactly $\CQ\inv(P)$ and the discussion in
Section~\ref{sec:cayley} shows that
\[
\iq([P])=\#(\CQ\inv_r(P)/\Stab(P))+2\,\rank(\Stab(P))-2
\]
\[
=\#\partial\Cayl(P)+2\,\rank(\Stab(P))-2=\sum_{P'\in [P]\cap
K_\CA}i_\CK(P').
\]

Adding up for all points $P$ in $K_\CA$, proves the Theorem.
\end{proof}

\subsection{Geometric index of a tree}

Gaboriau and Levitt in \cite{gl-rank} introduced the index of
$T$ as follows.

Let $P$ a point in $T$ and let $\pi_0(T\smallsetminus\{P\})$ be
the set of connected components of $T$ without $P$. The
stabilizer of $P$ acts on this set. The \textbf{geometric index}
of $P$ is
\[
\igeo(P)=\#(\pi_0(T\smallsetminus\{P\})/\Stab(P))+2\,\rank(\Stab(P))-2.
\]
This index is always non-negative because there are no
terminal vertices in a minimal tree. If the action of $\FN$ on
the tree $T$ is free the above definition becomes simpler:
\[
\igeo(P)=\#\pi_0(T\smallsetminus\{P\})-2.
\]
The geometric index
is constant inside an orbit under the action of $\FN$. The
\textbf{geometric index} of $T$ is then the sum of the indices
over all orbits of points:
\[
\igeo(T)=\sum_{[P]\in
T/\FN}\igeo(P).
\]

The following Theorem is proved by Gaboriau and Levitt:

\begin{thm}[\cite{gl-rank}]\label{thm:geometricindex} The
geometric index of an $\R$-tree with a very small minimal action
of the free group $\FN$ is bounded above by $2N-2$.\qed
\end{thm}

\subsection{Botanic of Trees}\label{sec:botanictree}

In this Section we establish a beginning of  classification of trees 
in the boundary of Outer Space.
Let $T$ be an $\R$-tree with a minimal very small action of $\FN$ by
isometries with dense orbits.
Let $\CA$ be a basis of $\FN$, let $\Omega_\CA$ be the relative
limit set, let $K_\CA$ be the compact heart of $T$ and let
$\CK_\CA=(K_\CA,\CA)$ be the associated system of isometries.

We first recall the existing terminology of geometric trees.
The tree $T$ is \textbf{geometric} if it can be obtained from
a system of isometries on a finite tree, as explained in the
Introduction (see for instance \cite{gab-indgen,best-survey}).
Geometric trees can be alternatively characterized thanks to the following:

\begin{thm}[\cite{gl-rank}, see also {\cite[Corollary~6.1]{chl4}}] 
\label{thm:geometric} 
Let $T$ be an $\R$-tree with a minimal, very small action by
isometries of $\FN$ with dense orbits. The following are equivalent:
\begin{enumerate}
\item $T$ is geometric;
\item the geometric index is maximal: $\igeo(T)=2N-2$;
\item $K_\CA$ is a finite tree. \qed
\end{enumerate}
\end{thm}

We now introduce more terminology.  The tree $T$ is of \textbf{surface
  type} if the Rips Machine, starting with the system of isometries
$\CK_\CA=(K_\CA,\CA)$, halts after finitely many steps.  More
precisely, a tree of surface type is:
\begin{itemize}
\item a \textbf{surface} tree if it is geometric
(this terminology is justified by the fact that a tree dual to measured
foliation on a surface with boundary is a surface tree.),
\item a \textbf{pseudo-surface} tree if it is not geometric.
\end{itemize}
(Note that the fact that a tree is pseudo-surface does not exactly 
correspond to the fact that the system of isometries $\CK_\CA=(K_\CA,\CA)$
is pseudo-surface, according to the definition given in 
Section~\ref{sec:pseudo-surface}).

By Theorem~\ref{thm:maxindex} and Theorem~\ref{thm:siiq} we get the
following characterization of trees of surface type:

\begin{thm}\label{thm:surface} Let $T$ be an $\R$-tree with a
minimal, very small action by isometries of $\FN$ with dense
orbits. The tree $T$ is of surface type if and only if its
$\CQ$-index is maximal: $\iq(T)=2N-2$. \qed
\end{thm}
This proves in particular that being of surface type is a property of $T$ 
and does not depend on the choice of a basis $\CA$ of the free group $\FN$.

The tree $T$ is  of \textbf{Levitt type}  if its relative
limit set $\Omega_\CA$ is totally disconnected (i.e. the
connected components of $\Omega_\CA$ are points). 
More precisely, a tree of Levitt type is
\begin{itemize}
\item \textbf{Levitt} if it is geometric (these trees were discovered by 
Levitt~\cite{levi-pseudogroup} and are also termed thin or exotic).
\item \textbf{pseudo-Levitt} if it is not geometric.
\end{itemize}
We now prove that being of Levitt type is a property of $T$ and does
not depend on the choice of a basis $\CA$ of $\FN$.

Let $L(T)$ be the dual lamination of the tree $T$. The
\textbf{limit set} $\Omega$ of $T$ is the image in the metric
completion $\bar T$ of $L(T)$ by the map $\CQ^2$:
\[
\Omega=\CQ^2(L(T)).
\]
Contrary to the relative limit set $\Omega_\CA$, the limit set
$\Omega$ is in general not closed.

\begin{thm}\label{thm:pseudolevittbasis} The tree $T$ is of
Levitt type if and only if the limit set $\Omega$ is totally
disconnected. 
\end{thm}
\begin{proof} 
  By definition, the relative limit set $\Omega_\CA=\CQ^2(L(T)\cap
  C_\CA(1))$ is a subset of $\Omega$. Thus, if $\Omega$ does not
  contain a non-trivial arc, $\Omega_\CA$ neither.

  Conversely, the double boundary of $\FN$ is the union of the
  translates of the unit cylinder and
\[
L(T)=\bigcup_{u\in\FN}u(L(T)\cap C_\CA(1))\ \text{ and }\
\Omega=\bigcup_{u\in\FN}u\Omega_\CA
\]
In particular if $I$ is a non-trivial arc in the limit set $\Omega$,
it is the countable union of its intersections with translates of the
relative limit set $\Omega_\CA$. Using Baire's Property for $I$, we
get that $\Omega_\CA$ contains a non-trivial arc.
\end{proof}

We remark that there are trees in the boundary of Outer Space
which are neither of surface or Levitt type. These are trees of
\textbf{mixed type}, that is to say their relative limit set
$\Omega_\CA$ contains non-trivial arcs but have infinitely many
connected components.

\begin{rem}
  A general classification of systems of isometries, in particular a
  Theorem à-la Imanishi to decompose trees of mixed type would be of
  interest.  More generally, the question of understanding the
  relationships between mixing properties of trees, indecomposability
  of systems of isometries and minimality of laminations seems to be
  natural. Together with Reynolds \cite{chr} we prove that
  indecomposable trees and minimal (up to diagonal leaves) laminations
  are dual to each other. In this spirit, see also the work of Reynolds~\cite{reyn}.
\end{rem}

\subsection{Mixing trees}\label{subsec:mixing}

In this section we give sufficient hypothesis on a tree to enforce
that it is either of surface type or of Levitt type. 

We first describe the limit set of trees of surface type.

\begin{prop}\label{prop:surfacetype} Let $T$ be an $\R$-tree
  with a minimal, very small action by isometries of $\FN$ with dense
  orbits. If the tree $T$ is of surface type, then the limit set
  $\Omega$ is connected and contains $T$.
\end{prop}
\begin{proof} As $T$ is minimal, and as $\Omega$ is $\FN$-invariant,
  we get that $\Omega$ is connected if and only if $\Omega$ contains
  $T$:
\[
T\subseteq\Omega\subseteq \bar T.
\]
If $T$ is of surface type, the Rips Machine starting with the system
of isometries $\CK_\CA=(K_\CA,\CA)$, halts after
finitely many steps, and thus $\Omega_\CA=F_n$ for some $n$, where
$F_n$ is the forest remaining after $n$ steps of the Rips Machine. The
system of isometries is reduced, hence the pseudo-orbit of each point
in $K_\CA$ is infinite and thus meets $F_n$, which proves that
$K_\CA\subseteq\FN.\Omega_\CA$. Moreover the orbit of each point in
$T$ meets $K_\CA$, thus $T\subseteq\FN K_\CA$, which concludes the
proof.
\end{proof}

A converse of this Proposition that requires stronger hypothesis
on $T$ is proved in Proposition~\ref{prop:mixingdichotomy}

A \textbf{segment} of an $\R$-tree is a subset isometric to a compact
real interval which is not reduced to a point.  The action of $\FN$ on
an $\R$-tree $T$ by isometries is \textbf{arc-dense} if every segment
of $T$ meets every orbits.  Following~\cite{morg-ergodic}, the action
of $\FN$ on an $\R$-tree $T$ by isometries is \textbf{mixing} if for
every segments $I$ and $J$ in $T$, the segment $J$ is covered by
finitely many translates of $I$: there exists finitely many elements
$u_1,\ldots,u_r$ of $\FN$ such that 
\[
J\subseteq u_1I\cup\cdots\cup u_rI. 
\]
It is obvious that a mixing action is arc-dense. An arc-dense
action has dense orbits and is minimal.

\begin{prop}\label{prop:mixingdichotomy} 
  Let $T$ be an $\R$-tree with a mixing action of $\FN$ by
  isometries. Then $T$ is either of surface type or of Levitt type.
\end{prop}
\begin{proof} 
  Let $\CA$ be a basis of $\FN$, let $\Omega_\CA$ be the relative
  limit set, $K_\CA$ the compact heart of $T$ and $S_\CA=(K_\CA,\CA)$
  be the associated system of isometries.

  By contradiction assume that $T$ is neither of surface type or of
  Levitt type. Then by definition of Levitt type, $\Omega_\CA$ contains a non-trivial
  connected component and thus a segment $I$. 

  Let $S_0=S_\CA$ and let $S_n=(F_n,A_n)$ be the sequence of systems
  of isometries obtained from $S_\CA$ by applying the Rips Machine. By
  definition of surface type, the Rips Machine runs forever.

Let $E_0$ be the set of points of $F_0=K_\CA$ erased at the
first step of the Rips Machine: $E_0=F_0\smallsetminus F_1$. As
$S_0$ is reduced, $E_0$ is contained in the convex hull in $F_0$
of $F_1$ and is a finite union of finite trees and an open
subset of $F_0$.

Let $E_n=F_n\smallsetminus F_{n+1}$ be the subset of $F_n$
erased at the $n+1$ step of the Rips Machine.
Let $D_n$ be the subset of $E_0$ defined by
\[
D_n=\{ P\in E_0\
|\ \exists u\in\FN, |u|=n\text{ and }P.u_i\in E_i,\text{ for
}i=1,\ldots,n\} 
\]
where $u_i$ is the prefix of $u$ of length $i$. By definition, for
each $n$, $D_{n+1}\subset D_n$. As the Rips Machine runs forever,
$D_n$ is a non-empty subset of $F_0$.
We distinguish two cases.

First assume that the nested intersection of the open non-empty
subsets $D_n$ is non-empty and let $P_0\in\cap_{n\in\N}D_n$. As
$P_0$ is in the open subset $E_0$ it is not an extremal point of
$F$ and as $T$ is arc-dense, there exists $u\in\FN$ such that
$uP_0\in I$. From Proposition~\ref{prop:chl4}, the partial
isometry $u\inv$ is defined at $P_0$ and $P_0u\inv=uP_0$ is in
the relative limit set $\Omega_\CA$. By definition of the Rips
Machine, for $n$ bigger than $|u|$, $P_0$ is not in $D_n$.
A contradiction.

Assume now that the nested intersection $\cap_{n\in\N}D_n$ is
empty and let $P_0$ be in the nested intersection of compact
subsets 
\[
P_0\in\bigcap_{n\in\N}\bar D_n. 
\]

Then, there exists $n_0$, such that for $n$ bigger than $n_0$, $P_0$
is in $\bar D_n\smallsetminus D_n$. Recall that $\bar D_{n_0}$ is a
finite tree and let $Q$ be a point of $D_{n_0}$. Then $Q$ is not an
extremal point of $K_\CA$ and let $J=[P_0,Q]$. The segment $J$ intersects
all the $D_n$ for $n\in\N$. If $P_0$ is in $T$ (and not in $\bar
T\smallsetminus T$) then, as $T$ is mixing, there there exist
$u_1,\ldots,u_r\in\FN$ such that
\[
J\subseteq u_1I\cup\cdots\cup u_rI.
\]
The partial isometries $u_1,\ldots,u_r$ are not empty and for each
$k$, $u_kI\cap K_\CA=Iu_k\inv$ thus, using
Proposition~\ref{prop:chl4},
\[
J\subseteq Iu_1\inv\cup\cdots\cup Iu_r\inv.
\]
Hence, for $n$ bigger than all the lengths of the $u_i$, $J\cap D_n$ is
empty.  A contradiction.

Thus, $P_0$ is in $\bar T\smallsetminus T$. We get that $P_0$ is
an extremal point of $\bar E_0$ and, as $E_0$ is open in
$K_\CA$, $P_0$ is in $\bar E_0\smallsetminus E_0$. As $P_0$ is
not in $E_0$ there are at least two partial isometries $a,b\in
A^{\pm 1}$ defined at $P_0$. One at least of $a$ and $b$ is not
defined in $E_0$, say $a$, and thus is defined only at $P_0$.
Thus, in $T$, there are at least two directions going out from
$P_0$ (one containing $K_\CA$ and the other containing $a\inv
K_\CA$). This contradicts the fact that $P_0$ is in $\bar
T\smallsetminus T$.
\end{proof}

\begin{cor}\label{cor:mixingsurface} Let $T$ be an $\R$-tree
with a mixing action of $\FN$ by isometries. Let $\Omega$ be the
limit set of $T$. The following are equivalent:
\begin{enumerate} 
\item\label{item:surf} $T$ is of surface type;
\item\label{item:omegaAforest}  $\Omega_\CA$ has finitely many components (i.e. $\Omega_\CA$ is a finite forest);
\item\label{item:conn} $\Omega$ is connected;
\item\label{item:T} $\Omega$ contains $T$, that is to say
$T\subseteq\Omega\subseteq\bar T$. 
\end{enumerate} 
\end{cor}
\begin{proof} 
  The equivalence of conditions~\ref{item:omegaAforest} and \ref{item:surf} 
  follows from the definition of
  surface type. Conditions~\ref{item:conn} and \ref{item:T} are
  equivalent because $T$ is minimal. We proved in
  Proposition~\ref{prop:surfacetype} that condition~\ref{item:surf}
  implies condition~\ref{item:conn}.  From the dichotomy of
  Proposition~\ref{prop:mixingdichotomy} and from
  Theorem~\ref{thm:pseudolevittbasis} we get that
  condition~\ref{item:conn} implies condition~\ref{item:surf}
\end{proof}

\section{Botanic memo}

In this Section we give a glossary of our classification of trees
for the working mathematician.

Let $T$ be an $\R$-tree with a minimal, very small action of
$\FN$ by isometries with dense orbits. We assume that the action
is indecomposable (see the definition in
Guirardel~\cite{guir-indecomposable}) or at least that $T$ is
not of mixed type.

For a basis $\CA$ of $\FN$, $\Omega_\CA$ is the relative limit
set and $K_\CA=\conv(\Omega_\CA)$ is the compact heart.
The compact heart $K_\CA$ is either a finite tree (and in the
good cases an interval) or not. This dichotomy is a property of
$T$ and does not depend of the choice of a particular basis
$\CA$ of $\FN$.

As we assumed that $T$ is indecomposable, $\Omega_\CA$ is either
a compact forest, that is to say it has finitely many connected
components (and in the good cases $\Omega_\CA=K_\CA$ is a tree)
or totally disconnected (and in the good cases a Cantor set).
This dichotomy is a property of $T$ and does not depend on the
choice of a particular basis $\CA$ of $\FN$.

The limit set of $T$ is $\Omega=\CQ^2(L(T))=\FN\cdot\Omega_\CA$. 
It is either totally disconnected (if $\Omega_\CA$ is), 
or it is connected (if $\Omega_\CA$ is a forest):
in the later case, $\Omega$ is  a tree,
$T\subseteq\Omega\subseteq\bar T$.

For such a tree $T$ we considered two indices: the geometric
index $\igeo(T)$ and the $\CQ$-index $\iq(T)$. Both indices
are bounded above by $2N-2$.
We sum up the terminology for $T$ and the results of
Section~\ref{sec:botanictree} in the following table.

\tabcolsep=.2em \noindent
\begin{tabular}{|c|c||c|c|} 
\cline{3-4}
\multicolumn{2}{c||}{} & \textbf{geometric} & \textbf{not geometric}\\ 
\cline{3-4}
\multicolumn{2}{c||}{} &
	\begin{tabular}{c}
	$K_\CA$ is a finite tree\\ 
	$\Updownarrow$\\
	$\igeo(T)=2N-2$
	\end{tabular} & 
		\begin{tabular}{c}
		$K_\CA$ is not a finite tree\\ 
		$\Updownarrow$\\ $\igeo(T)<2N-2$
		\end{tabular}\\
\hline\hline
\centering\rotatebox{90}
	{\begin{tabular}[b]{c}
	\hspace{-0.6cm}\textbf{Surface}\\
	\hspace{-0.6cm}\textbf{type}
	\end{tabular}} & 
		\begin{tabular}{c}
		$\Omega$ is a tree $(T\subset\Omega)$\\ 
		$\Updownarrow$\\
		$\Omega_\CA$ is a finite forest\\ 
		$\Updownarrow$\\
		$\iq(T)=2N-2$
		\end{tabular} &
			\textbf{surface} & \textbf{pseudo-surface}\\
\hline
\rotatebox{90}
	{\begin{tabular}{c}
	\hspace{-0.6cm}\textbf{Levitt}\\
	\hspace{-0.6cm}\textbf{type}
	\end{tabular}} & 
		\begin{tabular}{c}
		$\Omega$ is totally disconnected\\ 
		$\Updownarrow$\\
		$\Omega_\CA$ is totally disconnected\\ 
		$\Updownarrow$\\
		$\iq(T)<2N-2$
		\end{tabular} 
			& \textbf{Levitt} & \textbf{pseudo-Levitt}\\
\hline 
\end{tabular}

\end{document}